\documentclass[11pt]{amsart}
\usepackage{amsfonts,amscd,amssymb,amsmath,amsthm,mathrsfs,xcolor,lscape,amsmath,amssymb,latexsym}
\usepackage[utf8]{inputenc}
\usepackage[english]{babel}
\newtheorem{thm}{Theorem}

\newtheorem{rmk}{Remark}
\numberwithin{equation}{section}
\numberwithin{rmk}{section}

\newcommand{\C}{\mathbb{C}}

\newcommand{\Sp}{\mathrm{Sp}}

\newcommand{\GL}{\mathrm{GL}}

\newcommand{\Q}{\mathbb{Q}}
\newcommand{\Z}{\mathbb{Z}}
\newcommand{\e}{\epsilon}

\title{On Orthogonal Hypergeometric Groups of Degree Five}
\author{Jitendra Bajpai and Sandip Singh}
\address{Max Planck Institute for Mathematics, Vivatsgasse 7, 53111 Bonn, Germany\vspace{-1em}}
\curraddr{Mathematisches Institut, Georg-August Universit\"at G\"ottingen, Germany\vspace{-1em}}
\email{jitendra@math.uni-goettingen.de} 
\address{Department of Mathematics, Indian Institute of Technology Bombay, Mumbai, India\vspace{-1em}} \email{sandip@math.iitb.ac.in}

\begin{document}

\begin{abstract}
A computation shows that there are $77$ (up to scalar shifts) possible pairs of integer coefficient polynomials of degree five, having roots of unity as their roots, and satisfying the conditions of Beukers and Heckman \cite{BH89}, so that the Zariski closures of the associated monodromy groups are either finite or the orthogonal groups of non-degenerate and non-positive quadratic forms. Following the criterion of Beukers and Heckman \cite{BH89}, it is easy to check that only $4$ of these pairs correspond to finite monodromy groups and only $17$ pairs correspond to monodromy groups, for which, the Zariski closures have real rank one. 

There are remaining $56$ pairs, for which, the Zariski closures of the associated monodromy groups have real rank two. It follows from Venkataramana \cite{Ve14} that $11$ of these $56$ pairs correspond to arithmetic monodromy groups and the arithmeticity of $2$ other cases follows from Singh \cite{Si15}. In this article, we show that $23$ of the remaining $43$ rank two cases correspond to arithmetic groups.
\end{abstract}
\maketitle
%\setcounter{tocdepth}{2}
%\tableofcontents
%\vspace*{2em}

\section{Introduction}
We consider the monodromy action of the fundamental group $\pi_1$ of $\mathbb{P}^1(\C)\backslash\{0,1,\infty\}$ on the solution space of the $_n\mathrm{F}_{n-1}$ type hypergeometric differential equation 
\begin{equation}\label{introdifferentialequation}
D(\alpha;\beta)w=0
\end{equation}
on $\mathbb{P}^1(\C)$, having regular singularities at the points $\{0,1,\infty\}$, where the differential operator $D(\alpha;\beta)$ is defined as the following:
\begin{align*}
D(\alpha;\beta)&:=(\theta+\beta_1-1)\cdots(\theta+\beta_n-1)-z(\theta+\alpha_1)\cdots(\theta+\alpha_n)
\end{align*}
where $\theta=z\frac{d}{dz}$, and $\alpha=(\alpha_1,\alpha_2,\ldots,\alpha_n)$, $\beta=(\beta_1,\beta_2,\ldots,\beta_n)\in\Q^n$.

The monodromy groups of the hypergeometric differential equations (cf. Equation (\ref{introdifferentialequation})) are also called the hypergeometric groups, which are characterized by the following theorem of Levelt (\cite{Le61}; cf. \cite[Theorem 3.5]{BH89}):
\begin{thm}[Levelt]
 If $\alpha_1,\alpha_2,\ldots,\alpha_n$, $\beta_1,\beta_2,\ldots,\beta_n\in\Q$ such that $\alpha_j-\beta_k\not\in\Z$, for all $j,k=1,2,\ldots,n$, then there exists a basis of the solution space of the differential equation (\ref{introdifferentialequation}), with respect to which, the hypergeometric group is a subgroup of $\GL_n(\C)$ generated by the companion matrices $A$ and $B$ of  the polynomials 
\begin{equation}\label{parameterstopolynomials}
 f(X)=\prod_{j=1}^{n}(X-{\rm{e}^{2\pi i\alpha_j}})\quad\mbox{ and }\quad g(X)=\prod_{j=1}^{n}(X-{\rm{e}^{2\pi i\beta_j}})
\end{equation}
respectively, and the monodromy is defined by  $g_\infty\mapsto A$, $g_0\mapsto B^{-1}$, $g_1\mapsto A^{-1}B$, where $g_0, g_1, g_\infty$ are, respectively, the loops around $0,1,\infty$, which generate the $\pi_1$ of $\mathbb{P}^1(\C)\backslash\{0,1,\infty\}$ modulo the relation $g_\infty g_1 g_0=1$.
\end{thm}
We now denote the hypergeometric group by $\Gamma(f,g)$ which is a subgroup of $\GL_n(\C)$ generated by the companion matrices of the polynomials $f,g$. We consider the cases where the coefficients of $f,g$ are integers with $f(0)=\pm1$, $g(0)=\pm1$ (for example, one may take $f,g$ as product of cyclotomic polynomials); in these cases, $\Gamma(f,g)\subset\GL_n(\Z)$. We also assume that $f,g$ form a primitive pair \cite[Definition 5.1]{BH89}, are self-reciprocal and do not have any common root.
 
Beukers and Heckman \cite[Theorem 6.5]{BH89} have completely determined the Zariski closures $\mathrm{G}$ of the hypergeometric groups $\Gamma(f,g)$ which is briefly summarized as follows:

\begin{itemize}
\item If $n$ is even and $f(0)=g(0)=1$, then the hypergeometric group $\Gamma(f,g)$ preserves a non-degenerate integral symplectic form $\Omega$ on $\Z^n$ and $\Gamma(f,g)\subset\Sp_\Omega(\Z)$ is Zariski dense, that is, $\mathrm{G}=\Sp_\Omega$.

\item If $\Gamma(f,g)$ is {\it infinite} and $\frac{f(0)}{g(0)}=-1$, then $\Gamma(f,g)$ preserves a non-degenerate integral quadratic form $\mathrm{Q}$ on $\Z^n$ and $\Gamma(f,g)\subset\mathrm{O}_\mathrm{Q}(\Z)$ is Zariski dense, that is, $\mathrm{G}=\mathrm{O}_\mathrm{Q}$.

\item It follows from \cite[Corollary 4.7]{BH89} that $\Gamma(f,g)$ is {\it finite} if and only if either $\alpha_1<\beta_1<\alpha_2<\beta_2<\cdots<\alpha_n<\beta_n$ or $\beta_1<\alpha_1<\beta_2<\alpha_2<\cdots<\beta_n<\alpha_n$; and in this case, we say that the roots of $f$ and $g$ {\it interlace} on the unit circle.
\end{itemize}
  
 We call a hypergeometric group $\Gamma(f,g)$ {\it arithmetic}, if it is of {\it finite} index in $\mathrm{G}(\Z)$; and {\it thin}, if it has {\it infinite} index in $\mathrm{G}(\Z)$ \cite{Sa14}. 
 
 The following question of Sarnak \cite{Sa14} has drawn attention of many people: determine the pairs of polynomials $f,g$, for which, the associated hypergeometric group $\Gamma(f,g)$ is arithmetic. There have been some progress to answer this question.
 
 For the symplectic cases: infinitely many arithmetic $\Gamma(f,g)$ in $\Sp_n$ (for any even $n$) are given by Singh and Venkataramana in \cite{SiVe14}; some other examples of arithmetic $\Gamma(f,g)$ in $\Sp_4$ are given by Singh in \cite{Si14} and \cite{Si15S}; $7$ examples of thin $\Gamma(f,g)$ in $\Sp_4$ are given by Brav and Thomas in \cite{BT14}; and in \cite{HvS}, Hofmann and van Straten  have determined the index of $\Gamma(f,g)$ in $\Sp_4(\Z)$ for some of the arithmetic cases of \cite{Si14} and \cite{SiVe14}.
 
 For the orthogonal cases: when the quadratic form $\mathrm{Q}$ has signature $(n-1,1)$, Fuchs, Meiri and Sarnak give $7$ (infinite) families (depending on $n$ odd and $\geq5$) of thin $\Gamma(f,g)$ in \cite{FMS14}; when the quadratic form $\mathrm{Q}$ has signature $(p,q)$ with $p,q\geq2$, $11$ (infinite) families of arithmetic $\Gamma(f,g)$ are given by Venkataramana in \cite{Ve14}; an example of thin $\Gamma(f,g)$ in $\mathrm{O}(2,2)$ is given by Fuchs in \cite{Fu14}; and $2$ examples of arithmetic $\Gamma(f,g)$ in $\mathrm{O}(3,2)$ are given by Singh in \cite{Si15}, which deals with the $14$ orthogonal hypergeometric groups of degree five with a maximally unipotent monodromy, and these cases were inspired by the $14$ symplectic hypergeometric groups associated to Calabi-Yau threefolds \cite{Si14}.
 
 This article resulted from our effort to first write down all possible pairs $f,g$ (up to scalar shifts) of degree five integer coefficient polynomials, having roots of unity as their roots, forming a primitive pair and satisfying the condition: $f(0)=-1$ and $g(0)=1$; and then to answer the question to determine the arithmeticity or thinness of the associated hypergeometric groups $\Gamma(f,g)$. Note that, these pairs of the polynomials $f,g$ satisfy the conditions of Beukers and Heckman \cite{BH89} and hence the associated hypergeometric group $\Gamma(f,g)$ is either finite or preserves a non-degenerate integral quadratic form $\mathrm{Q}$ on $\Z^5$, and $\Gamma(f,g)\subset\mathrm{O}_\mathrm{Q}(\Z)$ is Zariski dense in the orthogonal group $\mathrm{O}_\mathrm{Q}$ of the quadratic form $\mathrm{Q}$.
 
 It is clear that there are finitely many pairs of degree five integer coefficient polynomials having roots of unity as their roots; among these pairs we find that there are $77$ pairs (cf. Tables \ref{table:realrank1}, \ref{table:arithmetic-Venkataramana}, \ref{table:arithmetic-Singh}, \ref{table:arithmetic-Bajpai-Singh}, \ref{table:finite}, \ref{table:unknown-Bajpai-Singh}, \ref{table:unknownrealrankone}) which satisfy the conditions of the last paragraph. Now, we consider the question to determine the {\it arithmeticity} or {\it thinness} of the associated orthogonal hypergeometric groups $\Gamma(f,g)$.
 
We note that the roots of the polynomials $f$ and $g$ interlace for the $4$ pairs $f,g$ of Table \ref{table:finite} and hence the hypergeometric groups associated to these pairs are finite \cite[Corollary 4.7]{BH89}. The quadratic forms associated to the $17$ pairs of Tables \ref{table:realrank1} and \ref{table:unknownrealrankone} have signature $(4,1)$, and the {\it thinness} of $\Gamma(f,g)$ associated to the $7$ pairs of Table \ref{table:realrank1} follows from Fuchs, Meiri, and Sarnak \cite{FMS14}. The arithmeticity or thinness of $\Gamma(f,g)$ associated to the other $10$ pairs of Table \ref{table:unknownrealrankone} is still {\it unknown}.

We also note that the quadratic forms associated to the remaining $56$ pairs (cf. Tables \ref{table:arithmetic-Venkataramana}, \ref{table:arithmetic-Singh}, \ref{table:arithmetic-Bajpai-Singh}, \ref{table:unknown-Bajpai-Singh}) have the signature $(3,2)$ and the arithmeticity of $\Gamma(f,g)$ associated to the $11$ pairs of Table \ref{table:arithmetic-Venkataramana} follows from Venkataramana \cite{Ve14}; and the arithmeticity of $\Gamma(f,g)$ associated to the other $2$ pairs of Table \ref{table:arithmetic-Singh} follows from Singh \cite{Si15}. Therefore, before this article arithmeticity of only $13$ hypergeometric groups $\Gamma(f,g)\subset\mathrm{O}(3,2)$ was known and the arithmeticity or thinness of the remaining $43$ hypergeometric groups $\Gamma(f,g)\subset\mathrm{O}(3,2)$ was {\it unknown}.

In this article, we show the {\it arithmeticity} of more than {\it half} of the remaining $43$ hypergeometric groups $\Gamma(f,g)\subset\mathrm{O}(3,2)$. In fact, we obtain the following theorem:

 \begin{thm}\label{maintheorem}
The hypergeometric groups associated to the $23$ pairs of polynomials in Table \ref{table:arithmetic-Bajpai-Singh} are arithmetic.
\end{thm}

Therefore, $23$ of the remaining $43$ pairs correspond to arithmetic hypergeometric groups; and the arithmeticity or thinness of the hypergeometric groups associated to the remaining $20$ pairs (cf. Table \ref{table:unknown-Bajpai-Singh}) is still {\it unknown}.
\vspace{-0.7em}
\section*{Acknowledgements}\vspace{-0.7em}
We thank Max Planck Institute for Mathematics for the postdoctoral fellowships, and for providing us very pleasant hospitalities. We thank the Mathematisches Forschungsinstitut Oberwolfach where we met T.N. Venkataramana, and thank him for many discussions on the subject of the paper.  We thank Wadim Zudilin for discussions at MPI. We also thank the referee for valuable comments and suggestions.

JB would like to extend his thanks to the Mathematisches Institut, Georg-August Universit\"at G\"ottingen for their support and hospitality during the final revision of the article. JB’s work was financially supported by ERC Consolidator grant 648329 (GRANT).

SS gratefully acknowledges the support of the INSPIRE Fellowship from the Department of Science and Technology, India. 
\vspace{-0.7em}

\section{Tables}\vspace{-0.7em}
In this section, we list all possible (up to scalar shifts) pairs $f,g$ of degree five polynomials, which are products of cyclotomic polynomials, and for which, the pair $(f,g)$ form a primitive pair \cite{BH89}, and $\frac{f(0)}{g(0)}=-1$, so that the associated monodromy group $\Gamma(f,g)$ preserves a quadratic form $\mathrm{Q}$, and $\Gamma(f,g)\subset\mathrm{O}_\mathrm{Q}(\Z)$ as a Zariski dense subgroup (except for the $4$ cases of Table \ref{table:finite}). Note that once we know the parameters $\alpha, \beta$, the polynomials $f, g$ can be determined by using the Equation (\ref{parameterstopolynomials}).

{\renewcommand{\arraystretch}{1.7}
\begin{table}[h]
\small
\addtolength{\tabcolsep}{-1pt}
\caption{List of the $7$ monodromy groups, for which, the associated quadratic form $\mathrm{Q}$ has signature $(4,1)$, and {\it thinness} follows from Fuchs, Meiri, and Sarnak \cite{FMS14}.}
\newcounter{rownum-0}
\setcounter{rownum-0}{0}
\centering
\begin{tabular}{ |c|  c|   c| c| c| c| c|}
\hline

  No. &$\alpha$ & $\beta$ &  No. & $\alpha$ &$\beta$\\ \hline
   \refstepcounter{rownum-0}\arabic{rownum-0}\label{thin1-rank1} & $\big( 0,0,0,\frac{1}{3},\frac{2}{3}\big)$ &$\big(\frac{1}{2},\frac{1}{8},\frac{3}{8},\frac{5}{8}, \frac{7}{8} \big)$& \refstepcounter{rownum-0}\arabic{rownum-0}\label{thin2-rank1} & $\big( 0,0,0,\frac{1}{3},\frac{2}{3}\big)$ &$\big(\frac{1}{2},\frac{1}{12},\frac{5}{12},\frac{7}{12}, \frac{11}{12} \big)$\\ \hline
   \refstepcounter{rownum-0}\arabic{rownum-0}\label{thin3-rank1} & $\big(0,0,0,\frac{1}{4},\frac{3}{4}\big)$  & $\big(\frac{1}{2},\frac{1}{10},\frac{3}{10},\frac{7}{10}, \frac{9}{10} \big)$&\refstepcounter{rownum-0}\arabic{rownum-0}\label{thin4-rank1} & $\big(0,0,0,\frac{1}{6},\frac{5}{6}\big)$  & $\big(\frac{1}{2},\frac{1}{10},\frac{3}{10},\frac{7}{10}, \frac{9}{10} \big)$\\ \hline
   \refstepcounter{rownum-0}\arabic{rownum-0}\label{thin5-rank1} & $\big(0,\frac{1}{3},\frac{2}{3},\frac{1}{6},\frac{5}{6}\big)$ & $\big(\frac{1}{2},\frac{1}{5},\frac{2}{5},\frac{3}{5},\frac{4}{5} \big)$&\refstepcounter{rownum-0}\arabic{rownum-0}\label{thin6-rank1} & $\big(0,\frac{1}{10},\frac{3}{10},\frac{7}{10},\frac{9}{10}\big)$  & $\big(\frac{1}{2},\frac{1}{3},\frac{2}{3},\frac{1}{4},\frac{3}{4}\big)$\\ \hline
   \refstepcounter{rownum-0}\arabic{rownum-0}\label{thin7-rank1} & $\big(0,\frac{1}{10},\frac{3}{10},\frac{7}{10},\frac{9}{10}\big)$  & $\big(\frac{1}{2},\frac{1}{8},\frac{3}{8},\frac{5}{8}, \frac{7}{8} \big)$ & & &\\ \hline
   
\end{tabular}
\label{table:realrank1}
\end{table}}
%\vspace{-0cm}    

{\renewcommand{\arraystretch}{1.7}
\begin{table}[h]
\small
\addtolength{\tabcolsep}{-1pt}
\caption{List of the $11$ monodromy groups, for which, the associated quadratic form $\mathrm{Q}$ has signature $(3,2)$, and {\it arithmeticity} follows from Venkataramana \cite{Ve14}.}
\centering
\begin{tabular}{ |c|  c|   c| c| c| c|}
\hline

  No. &$\alpha$ & $\beta$ & No. &$\alpha$ & $\beta$\\ \hline
   \refstepcounter{rownum-0}\arabic{rownum-0}\label{arithmetic-Venkataramana1} & $\big( 0,0,0,\frac{1}{3},\frac{2}{3}\big)$  &$\big(\frac{1}{2},\frac{1}{10},\frac{3}{10},\frac{7}{10}, \frac{9}{10} \big)$&\refstepcounter{rownum-0}\arabic{rownum-0}\label{arithmetic-Venkataramana2} & $\big(0,\frac{1}{3},\frac{1}{3},\frac{2}{3},\frac{2}{3}\big)$  &$\big(\frac{1}{2},\frac{1}{4},\frac{1}{4},\frac{3}{4},\frac{3}{4}\big)$\\ \hline
     \refstepcounter{rownum-0}\arabic{rownum-0}\label{arithmetic-Venkataramana3} & $\big(0,\frac{1}{3},\frac{1}{3},\frac{2}{3},\frac{2}{3}\big)$  & $\big(\frac{1}{2},\frac{1}{5},\frac{2}{5},\frac{3}{5},\frac{4}{5} \big)$&\refstepcounter{rownum-0}\arabic{rownum-0}\label{arithmetic-Venkataramana4} & $\big(0,\frac{1}{3},\frac{1}{3},\frac{2}{3},\frac{2}{3}\big)$  & $\big(\frac{1}{2}, \frac{1}{6},\frac{1}{6},\frac{5}{6},\frac{5}{6} \big)$\\ \hline
      \refstepcounter{rownum-0}\arabic{rownum-0}\label{arithmetic-Venkataramana5} & $\big(0,\frac{1}{3},\frac{1}{3},\frac{2}{3},\frac{2}{3}\big)$  &$\big(\frac{1}{2},\frac{1}{10},\frac{3}{10},\frac{7}{10}, \frac{9}{10} \big)$&\refstepcounter{rownum-0}\arabic{rownum-0}\label{arithmetic-Venkataramana6} & $\big(0,\frac{1}{3},\frac{1}{3},\frac{2}{3},\frac{2}{3}\big)$  &$\big(\frac{1}{2},\frac{1}{12},\frac{5}{12},\frac{7}{12}, \frac{11}{12} \big)$\\ \hline
      \refstepcounter{rownum-0}\arabic{rownum-0}\label{arithmetic-Venkataramana7} & $\big(0,\frac{1}{3},\frac{2}{3},\frac{1}{4},\frac{3}{4}\big)$   &$\big(\frac{1}{2},\frac{1}{5},\frac{2}{5},\frac{3}{5},\frac{4}{5} \big)$&\refstepcounter{rownum-0}\arabic{rownum-0}\label{arithmetic-Venkataramana8} & $\big(0,\frac{1}{5},\frac{2}{5},\frac{3}{5},\frac{4}{5}\big)$  & $\big(\frac{1}{2},\frac{1}{4},\frac{1}{4},\frac{3}{4},\frac{3}{4}\big)$\\ \hline
 \refstepcounter{rownum-0}\arabic{rownum-0}\label{arithmetic-Venkataramana9} &$\big(0,\frac{1}{10},\frac{3}{10},\frac{7}{10},\frac{9}{10}\big)$&$\big(\frac{1}{2}, \frac{1}{4},\frac{1}{4},\frac{3}{4},\frac{3}{4} \big)$&\refstepcounter{rownum-0}\arabic{rownum-0}\label{arithmetic-Venkataramana10} &$\big(0,\frac{1}{12},\frac{5}{12},\frac{7}{12},\frac{11}{12}\big)$ & $\big(\frac{1}{2}, \frac{1}{4},\frac{1}{4},\frac{3}{4},\frac{3}{4} \big)$\\ \hline
\refstepcounter{rownum-0}\arabic{rownum-0}\label{arithmetic-Venkataramana11} &$\big(0,\frac{1}{12},\frac{5}{12},\frac{7}{12},\frac{11}{12}\big)$ & $\big(\frac{1}{2},\frac{1}{10},\frac{3}{10},\frac{7}{10}, \frac{9}{10} \big)$& & &\\ \hline
 \end{tabular}
\label{table:arithmetic-Venkataramana}
\end{table}}
\vspace{0.5cm}

{\renewcommand{\arraystretch}{1.7}
\begin{table}[h]
\small
\addtolength{\tabcolsep}{-1pt}
\caption{List of the $2$ monodromy groups, for which, the associated quadratic form $\mathrm{Q}$ has signature $(3,2)$, and {\it arithmeticity} follows from Singh \cite{Si15}.}
\centering
\begin{tabular}{ |c|  c|   c| c| c| c|}
\hline

  No. &$\alpha$ & $\beta$ & No. &$\alpha$ & $\beta$\\ \hline
   \refstepcounter{rownum-0}\arabic{rownum-0}\label{arithmetic-Singh1} & $\big(0,0,0,0,0\big)$  &$\big(\frac{1}{2},\frac{1}{4},\frac{1}{4},\frac{3}{4},\frac{3}{4}\big)$&\refstepcounter{rownum-0}\arabic{rownum-0}\label{arithmetic-Singh2} & $\big(0,0,0,0,0\big)$  & $\big(\frac{1}{2}, \frac{1}{6},\frac{1}{6},\frac{5}{6},\frac{5}{6} \big)$\\ \hline 
   
   \end{tabular}
\label{table:arithmetic-Singh}
\end{table}}

{\renewcommand{\arraystretch}{1.7}
\begin{table}[h]
\small
\addtolength{\tabcolsep}{-1pt}
\caption{List of the $23$ monodromy groups, for which, the associated quadratic form $\mathrm{Q}$ has signature $(3,2)$, and {\it arithmeticity} is shown in Section \ref{proof} of this paper.}
\centering
\begin{tabular}{ |c|  c|   c| c| c| c|}
\hline

  No. &$\alpha$ & $\beta$ & No. &$\alpha$ & $\beta$\\ \hline
  \refstepcounter{rownum-0}\arabic{rownum-0}\label{arithmetic-BS1} & $\big( 0,0,0,\frac{1}{3},\frac{2}{3}\big)$  &$\big(\frac{1}{2},\frac{1}{4},\frac{1}{4},\frac{3}{4},\frac{3}{4}\big)$ &\refstepcounter{rownum-0}\arabic{rownum-0}\label{arithmetic-BS2} & $\big( 0,0,0,\frac{1}{3},\frac{2}{3}\big)$  &$\big(\frac{1}{2},\frac{1}{4},\frac{3}{4},\frac{1}{6},\frac{5}{6}\big)$\\ \hline
  \refstepcounter{rownum-0}\arabic{rownum-0}\label{arithmetic-BS3} & $\big( 0,0,0,\frac{1}{3},\frac{2}{3}\big)$  &$\big(\frac{1}{2}, \frac{1}{6},\frac{1}{6},\frac{5}{6},\frac{5}{6} \big)$  &\refstepcounter{rownum-0}\arabic{rownum-0}\label{arithmetic-BS4} & $\big(0,0,0,\frac{1}{4},\frac{3}{4}\big)$   &$\big(\frac{1}{2}, \frac{1}{6},\frac{1}{6},\frac{5}{6},\frac{5}{6} \big)$\\ \hline
  \refstepcounter{rownum-0}\arabic{rownum-0}\label{arithmetic-BS5} & $\big(0,0,0,\frac{1}{6},\frac{5}{6}\big)$  &$\big(\frac{1}{2}, \frac{1}{4},\frac{1}{4},\frac{3}{4},\frac{3}{4} \big)$  &\refstepcounter{rownum-0}\arabic{rownum-0}\label{arithmetic-BS6} & $\big(0,\frac{1}{3},\frac{1}{3},\frac{2}{3},\frac{2}{3}\big)$ &$\big(\frac{1}{2},\frac{1}{4},\frac{3}{4},\frac{1}{6},\frac{5}{6}\big)$\\ \hline
  \refstepcounter{rownum-0}\arabic{rownum-0}\label{arithmetic-BS7} & $\big(0,\frac{1}{3},\frac{1}{3},\frac{2}{3},\frac{2}{3}\big)$ &$\big(\frac{1}{2},\frac{1}{8},\frac{3}{8},\frac{5}{8}, \frac{7}{8} \big)$ &\refstepcounter{rownum-0}\arabic{rownum-0}\label{arithmetic-BS8} & $\big(0,\frac{1}{3},\frac{2}{3},\frac{1}{4},\frac{3}{4}\big)$ &$\big(\frac{1}{2},\frac{1}{8},\frac{3}{8},\frac{5}{8}, \frac{7}{8} \big)$\\ \hline
   \refstepcounter{rownum-0}\arabic{rownum-0}\label{arithmetic-BS9} & $\big(0,\frac{1}{3},\frac{2}{3},\frac{1}{4},\frac{3}{4}\big)$ &$\big(\frac{1}{2},\frac{1}{12},\frac{5}{12},\frac{7}{12}, \frac{11}{12} \big)$ &\refstepcounter{rownum-0}\arabic{rownum-0}\label{arithmetic-BS10} & $\big(0,\frac{1}{3},\frac{2}{3},\frac{1}{6},\frac{5}{6}\big)$&$\big(\frac{1}{2}, \frac{1}{4},\frac{1}{4},\frac{3}{4},\frac{3}{4} \big)$\\ \hline
   \refstepcounter{rownum-0}\arabic{rownum-0}\label{arithmetic-BS11} & $\big(0,\frac{1}{3},\frac{2}{3},\frac{1}{6},\frac{5}{6}\big)$&$\big(\frac{1}{2},\frac{1}{8},\frac{3}{8},\frac{5}{8}, \frac{7}{8} \big)$ &\refstepcounter{rownum-0}\arabic{rownum-0}\label{arithmetic-BS12} & $\big(0,\frac{1}{3},\frac{2}{3},\frac{1}{6},\frac{5}{6}\big)$&$\big(\frac{1}{2},\frac{1}{12},\frac{5}{12},\frac{7}{12}, \frac{11}{12} \big)$\\ \hline
   \refstepcounter{rownum-0}\arabic{rownum-0}\label{arithmetic-BS13} & $\big(0,\frac{1}{5},\frac{2}{5},\frac{3}{5},\frac{4}{5}\big)$&$\big(\frac{1}{2},\frac{1}{3},\frac{1}{3},\frac{2}{3},\frac{2}{3}\big)$  &\refstepcounter{rownum-0}\arabic{rownum-0}\label{arithmetic-BS14} & $\big(0,\frac{1}{5},\frac{2}{5},\frac{3}{5},\frac{4}{5}\big)$&$\big(\frac{1}{2},\frac{1}{3},\frac{2}{3},\frac{1}{4},\frac{3}{4}\big)$\\ \hline
  \refstepcounter{rownum-0}\arabic{rownum-0}\label{arithmetic-BS15} & $\big(0,\frac{1}{6},\frac{1}{6},\frac{5}{6},\frac{5}{6}\big)$&$\big(\frac{1}{2},\frac{1}{3},\frac{1}{3},\frac{2}{3},\frac{2}{3}\big)$ &\refstepcounter{rownum-0}\arabic{rownum-0}\label{arithmetic-BS16} & $\big(0,\frac{1}{6},\frac{1}{6},\frac{5}{6},\frac{5}{6}\big)$&$\big(\frac{1}{2},\frac{1}{3},\frac{2}{3},\frac{1}{4},\frac{3}{4}\big)$\\ \hline
  \refstepcounter{rownum-0}\arabic{rownum-0}\label{arithmetic-BS17} & $\big(0,\frac{1}{6},\frac{1}{6},\frac{5}{6},\frac{5}{6}\big)$&$\big(\frac{1}{2},\frac{1}{4},\frac{1}{4}, \frac{3}{4},\frac{3}{4}\big)$ &\refstepcounter{rownum-0}\arabic{rownum-0}\label{arithmetic-BS18} & $\big(0,\frac{1}{8},\frac{3}{8},\frac{5}{8},\frac{7}{8}\big)$&$\big(\frac{1}{2},\frac{1}{3},\frac{1}{3},\frac{2}{3},\frac{2}{3}\big)$\\ \hline
  \refstepcounter{rownum-0}\arabic{rownum-0}\label{arithmetic-BS19} & $\big(0,\frac{1}{8},\frac{3}{8},\frac{5}{8},\frac{7}{8}\big)$&$\big(\frac{1}{2},\frac{1}{3},\frac{2}{3},\frac{1}{4},\frac{3}{4}\big)$&\refstepcounter{rownum-0}\arabic{rownum-0}\label{arithmetic-BS20} & $\big(0,\frac{1}{8},\frac{3}{8},\frac{5}{8},\frac{7}{8}\big)$&$\big(\frac{1}{2},\frac{1}{4},\frac{1}{4},\frac{3}{4},\frac{3}{4}\big)$\\ \hline
  \refstepcounter{rownum-0}\arabic{rownum-0}\label{arithmetic-BS21} & $\big(0,\frac{1}{12},\frac{5}{12},\frac{7}{12},\frac{11}{12}\big)$ &$\big(\frac{1}{2},\frac{1}{3},\frac{2}{3},\frac{1}{4},\frac{3}{4}\big)$ &\refstepcounter{rownum-0}\arabic{rownum-0}\label{arithmetic-BS22} & $\big(0,\frac{1}{12},\frac{5}{12},\frac{7}{12},\frac{11}{12}\big)$ &$\big(\frac{1}{2},\frac{1}{5},\frac{2}{5},\frac{3}{5},\frac{4}{5} \big)$\\ \hline
   \refstepcounter{rownum-0}\arabic{rownum-0}\label{arithmetic-BS23} & $\big(0,\frac{1}{12},\frac{5}{12},\frac{7}{12},\frac{11}{12}\big)$ &$\big(\frac{1}{2},\frac{1}{8},\frac{3}{8},\frac{5}{8}, \frac{7}{8} \big)$& & &\\ \hline
     \end{tabular}
\label{table:arithmetic-Bajpai-Singh}
\end{table}}
\vspace{2cm}

{\renewcommand{\arraystretch}{1.7}
\begin{table}[h]
\small
\addtolength{\tabcolsep}{-1pt}
\caption{List of the $4$ monodromy groups, for which, the associated quadratic form $\mathrm{Q}$ is positive definite (since the roots of the corresponding polynomials interlace on the unit circle \cite[Corollary 4.7]{BH89}).}
\centering
\begin{tabular}{ |c|  c|   c| c| c| c|}
\hline

  No. &$\alpha$ & $\beta$ & No.& $\alpha$ & $\beta$\\ \hline
   \refstepcounter{rownum-0}\arabic{rownum-0}\label{finite1-rank0} &$\big(0,\frac{1}{3},\frac{2}{3},\frac{1}{4},\frac{3}{4}\big)$ & $\big(\frac{1}{2},\frac{1}{10},\frac{3}{10},\frac{7}{10}, \frac{9}{10} \big)$&\refstepcounter{rownum-0}\arabic{rownum-0}\label{finite2-rank0} &$\big(0,\frac{1}{3},\frac{2}{3},\frac{1}{6},\frac{5}{6}\big)$ &$\big(\frac{1}{2},\frac{1}{10},\frac{3}{10},\frac{7}{10}, \frac{9}{10} \big)$\\ \hline 
   \refstepcounter{rownum-0}\arabic{rownum-0}\label{finite3-rank0} &$\big(0,\frac{1}{5},\frac{2}{5},\frac{3}{5},\frac{4}{5}\big)$  &$\big(\frac{1}{2},\frac{1}{10},\frac{3}{10},\frac{7}{10}, \frac{9}{10} \big)$&\refstepcounter{rownum-0}\arabic{rownum-0}\label{finite4-rank0} &$\big(0,\frac{1}{8},\frac{3}{8},\frac{5}{8},\frac{7}{8}\big)$ &$\big(\frac{1}{2},\frac{1}{10},\frac{3}{10},\frac{7}{10}, \frac{9}{10} \big)$\\ \hline 
 
  \end{tabular}
\label{table:finite}
\end{table}}
 \clearpage 
 
{\renewcommand{\arraystretch}{1.7}
\begin{table}[h]
\small
\addtolength{\tabcolsep}{-1pt}
\caption{List of the $20$ monodromy groups, for which, the associated quadratic form $\mathrm{Q}$ has signature $(3,2)$, and {\it arithmeticity} or {\it thinness} is {\it unknown}.}
\centering
\begin{tabular}{ |c|  c|   c| c| c| c|}
\hline
  No. &$\alpha$ & $\beta$ & No. &$\alpha$ & $\beta$\\ \hline
  
  \refstepcounter{rownum-0}\arabic{rownum-0}\label{unknown-BS1} & $\big( 0,0,0,0,0\big)$  & $\big(\frac{1}{2},\frac{1}{2},\frac{1}{2},\frac{1}{2},\frac{1}{2}\big)$ &\refstepcounter{rownum-0}\arabic{rownum-0}\label{unknown-BS2} & $\big( 0,0,0,0,0\big)$  & $\big(\frac{1}{2},\frac{1}{2},\frac{1}{2},\frac{1}{3},\frac{2}{3}\big)$\\ \hline
  \refstepcounter{rownum-0}\arabic{rownum-0}\label{unknown-BS3} &  $\big( 0,0,0,0,0\big)$ &  $\big(\frac{1}{2},\frac{1}{2},\frac{1}{2},\frac{1}{4},\frac{3}{4}\big)$&\refstepcounter{rownum-0}\arabic{rownum-0}\label{unknown-BS4} &$\big( 0,0,0,0,0\big)$ & $\big(\frac{1}{2},\frac{1}{2},\frac{1}{2},\frac{1}{6},\frac{5}{6}\big)$\\ \hline
  \refstepcounter{rownum-0}\arabic{rownum-0}\label{unknown-BS5} & $\big( 0,0,0,0,0\big)$& $\big(\frac{1}{2},\frac{1}{3},\frac{1}{3},\frac{2}{3},\frac{2}{3}\big)$ &\refstepcounter{rownum-0}\arabic{rownum-0}\label{unknown-BS6} & $\big( 0,0,0,0,0\big)$  &$\big(\frac{1}{2},\frac{1}{3},\frac{2}{3},\frac{1}{4},\frac{3}{4}\big)$\\ \hline 
  \refstepcounter{rownum-0}\arabic{rownum-0}\label{unknown-BS7} & $\big( 0,0,0,0,0\big)$   &$\big(\frac{1}{2},\frac{1}{3},\frac{2}{3},\frac{1}{6},\frac{5}{6}\big)$  &\refstepcounter{rownum-0}\arabic{rownum-0}\label{unknown-BS8} & $\big( 0,0,0,0,0\big)$   & $\big(\frac{1}{2},\frac{1}{4},\frac{3}{4},\frac{1}{6},\frac{5}{6}\big)$\\ \hline 
  \refstepcounter{rownum-0}\arabic{rownum-0}\label{unknown-BS9} &  $\big( 0,0,0,0,0\big)$  & $\big(\frac{1}{2},\frac{1}{5},\frac{2}{5},\frac{3}{5},\frac{4}{5} \big)$ &\refstepcounter{rownum-0}\arabic{rownum-0}\label{unknown-BS10} &$\big( 0,0,0,0,0\big)$ & $\big(\frac{1}{2},\frac{1}{8},\frac{3}{8},\frac{5}{8}, \frac{7}{8} \big)$\\ \hline
   \refstepcounter{rownum-0}\arabic{rownum-0}\label{unknown-BS11} & $\big( 0,0,0,0,0\big)$  &$\big(\frac{1}{2},\frac{1}{10},\frac{3}{10},\frac{7}{10}, \frac{9}{10} \big)$   &\refstepcounter{rownum-0}\arabic{rownum-0}\label{unknown-BS12} & $\big( 0,0,0,0,0\big)$ & $\big(\frac{1}{2},\frac{1}{12},\frac{5}{12},\frac{7}{12}, \frac{11}{12} \big)$\\ \hline
   \refstepcounter{rownum-0}\arabic{rownum-0}\label{unknown-BS14} & $\big(0,0,0,\frac{1}{4},\frac{3}{4}\big)$  & $\big(\frac{1}{2},\frac{1}{2},\frac{1}{2},\frac{1}{3},\frac{2}{3}\big)$&\refstepcounter{rownum-0}\arabic{rownum-0}\label{unknown-BS15} &$\big(0,0,0,\frac{1}{4},\frac{3}{4}\big)$ &$\big(\frac{1}{2},\frac{1}{3},\frac{1}{3},\frac{2}{3},\frac{2}{3}\big)$\\ \hline
   \refstepcounter{rownum-0}\arabic{rownum-0}\label{unknown-BS17} &$\big(0,0,0,\frac{1}{6},\frac{5}{6}\big)$ & $\big(\frac{1}{2},\frac{1}{2},\frac{1}{2},\frac{1}{3},\frac{2}{3}\big)$&\refstepcounter{rownum-0}\arabic{rownum-0}\label{unknown-BS18} &$\big(0,0,0,\frac{1}{6},\frac{5}{6}\big)$ & $\big(\frac{1}{2},\frac{1}{3},\frac{1}{3},\frac{2}{3},\frac{2}{3}\big)$\\ \hline
   \refstepcounter{rownum-0}\arabic{rownum-0}\label{unknown-BS19} &$\big(0,0,0,\frac{1}{6},\frac{5}{6}\big)$ & $\big(\frac{1}{2},\frac{1}{3},\frac{2}{3},\frac{1}{4},\frac{3}{4}\big)$&\refstepcounter{rownum-0}\arabic{rownum-0}\label{unknown-BS20} & $\big(0,0,0,\frac{1}{6},\frac{5}{6}\big)$  & $\big(\frac{1}{2},\frac{1}{5},\frac{2}{5},\frac{3}{5},\frac{4}{5} \big)$\\ \hline
   \refstepcounter{rownum-0}\arabic{rownum-0}\label{unknown-BS21} &$\big(0,\frac{1}{10},\frac{3}{10},\frac{7}{10},\frac{9}{10}\big)$ & $\big(\frac{1}{2},\frac{1}{3},\frac{1}{3},\frac{2}{3},\frac{2}{3}\big)$&\refstepcounter{rownum-0}\arabic{rownum-0}\label{unknown-BS23} &$\big(0,\frac{1}{12},\frac{5}{12},\frac{7}{12},\frac{11}{12}\big)$ & $\big(\frac{1}{2},\frac{1}{3},\frac{1}{3},\frac{2}{3},\frac{2}{3}\big)$\\ \hline
   \end{tabular}
\label{table:unknown-Bajpai-Singh}
\end{table}}
 \vspace{2cm}
{\renewcommand{\arraystretch}{1.7}
\begin{table}[h]
\small
\addtolength{\tabcolsep}{-1pt}
\caption{List of the $10$ monodromy groups, for which, the associated quadratic form $\mathrm{Q}$ has signature $(4,1)$, and {\it arithmeticity} or {\it thinness} is {\it unknown}.}
\centering
\begin{tabular}{ |c|  c|   c| c| c| c|}
\hline
  No. &$\alpha$ & $\beta$ & No. &$\alpha$ & $\beta$\\ \hline
  \refstepcounter{rownum-0}\arabic{rownum-0}\label{unknown1-rank1} & $\big( 0,0,0,\frac{1}{3},\frac{2}{3}\big)$ & $\big(\frac{1}{2},\frac{1}{2},\frac{1}{2},\frac{1}{4},\frac{3}{4}\big)$  &\refstepcounter{rownum-0}\arabic{rownum-0}\label{unknown-BS13} & $\big( 0,0,0,\frac{1}{3},\frac{2}{3}\big)$  & $\big(\frac{1}{2},\frac{1}{2},\frac{1}{2},\frac{1}{6},\frac{5}{6}\big)$\\ \hline
  \refstepcounter{rownum-0}\arabic{rownum-0}\label{unknown2-rank1} & $\big( 0,0,0,\frac{1}{3},\frac{2}{3}\big)$ &$\big(\frac{1}{2},\frac{1}{5},\frac{2}{5},\frac{3}{5},\frac{4}{5} \big)$&\refstepcounter{rownum-0}\arabic{rownum-0}\label{unknown3-rank1} &$\big(0,0,0,\frac{1}{4},\frac{3}{4}\big)$ &$\big(\frac{1}{2},\frac{1}{3},\frac{2}{3},\frac{1}{6},\frac{5}{6}\big)$\\ \hline
  \refstepcounter{rownum-0}\arabic{rownum-0}\label{unknown4-rank1} &$\big(0,0,0,\frac{1}{4},\frac{3}{4}\big)$ & $\big(\frac{1}{2},\frac{1}{5},\frac{2}{5},\frac{3}{5},\frac{4}{5} \big)$&\refstepcounter{rownum-0}\arabic{rownum-0}\label{unknown-BS16} &$\big(0,0,0,\frac{1}{4},\frac{3}{4}\big)$ & $\big(\frac{1}{2},\frac{1}{8},\frac{3}{8},\frac{5}{8}, \frac{7}{8} \big)$\\ \hline
  \refstepcounter{rownum-0}\arabic{rownum-0}\label{unknown5-rank1} &$\big(0,0,0,\frac{1}{4},\frac{3}{4}\big)$ & $\big(\frac{1}{2},\frac{1}{12},\frac{5}{12},\frac{7}{12}, \frac{11}{12} \big)$&\refstepcounter{rownum-0}\arabic{rownum-0}\label{unknown6-rank1} &$\big(0,0,0,\frac{1}{6},\frac{5}{6}\big)$  &$\big(\frac{1}{2},\frac{1}{8},\frac{3}{8},\frac{5}{8}, \frac{7}{8} \big)$\\ \hline
  \refstepcounter{rownum-0}\arabic{rownum-0}\label{unknown7-rank1} &$\big(0,0,0,\frac{1}{6},\frac{5}{6}\big)$  & $\big(\frac{1}{2},\frac{1}{12},\frac{5}{12},\frac{7}{12}, \frac{11}{12} \big)$&\refstepcounter{rownum-0}\arabic{rownum-0}\label{unknown-BS22} & $\big(0,\frac{1}{10},\frac{3}{10},\frac{7}{10},\frac{9}{10}\big)$& $\big(\frac{1}{2},\frac{1}{5},\frac{2}{5},\frac{3}{5},\frac{4}{5} \big)$\\ \hline
  \end{tabular}
\label{table:unknownrealrankone}
\end{table}}
\clearpage

\section{Proof of Theorem \ref{maintheorem}}\label{proof}
In this section, we show the arithmeticity of all the hypergeometric groups associated to the pairs of polynomials listed in Table \ref{table:arithmetic-Bajpai-Singh}, and it proves Theorem \ref{maintheorem}.
\subsection*{Strategy} We first compute the quadratic forms $\mathrm{Q}$ (up to scalar multiplications) preserved by the hypergeometric groups of Theorem \ref{maintheorem}, and then show that the real rank of the orthogonal group $\mathrm{SO}_\mathrm{Q}$ is two, and the $\Q$ - rank is either one or two. We then form a basis $\{\e_1,\e_2,\e_3,\e_2^*,\e_1^*\}$ of $\Q^5$, which satisfy the following condition: in $\Q$ - rank two cases, the matrix form of the quadratic form $\mathrm{Q}$, with respect to the basis $\{\e_1,\e_2,\e_3,\e_2^*,\e_1^*\}$, is anti-diagonal; and in $\Q$ - rank one cases, the vectors $\e_1, \e_1^*$ are $\mathrm{Q}$ - isotropic non-orthogonal vectors (that is, $\mathrm{Q}(\e_1,\e_1)=\mathrm{Q}(\e_1^*,\e_1^*)=0$ and $\mathrm{Q}(\e_1,\e_1^*)\neq0$), and  the vectors $\e_2,\e_3,\e_2^*$ are $\mathrm{Q}$ - orthogonal to the vectors $\e_1, \e_1^*$.

Let $\mathrm{P}$ be the parabolic $\Q$ - subgroup of $\mathrm{SO}_\mathrm{Q}$, which preserves the following flag:

\[\{0\}\subset\Q\e_1\subset\Q\e_1\oplus\Q\e_2\oplus\Q\e_3\oplus\Q\e_2^*\subset\Q\e_1\oplus\Q\e_2\oplus\Q\e_3\oplus\Q\e_2^*\oplus\Q\e_1^*\]\\
and $\mathrm{U}$ be the unipotent radical of $\mathrm{P}$. It can be checked easily that the unipotent radical $\mathrm{U}$ is isomorphic to $\Q^3$ (as a group), and in particular, $\mathrm{U}(\Z)$ is a free abelian group isomorphic to $\Z^3$. 

We prove Theorem \ref{maintheorem} (except for the cases \ref{arithmetic-BS16}, \ref{arithmetic-BS19}, and \ref{arithmetic-BS22} of Table \ref{table:arithmetic-Bajpai-Singh})  by using the following criterion of Raghunathan \cite{Rag91}, and Venkataramana \cite{Ve94} (cf. \cite[Theorem 6]{Ve14}): If $\Gamma(f,g)\subset\mathrm{O}_\mathrm{Q}(\Z)$ is Zariski dense and intersecting $\mathrm{U}(\Z)$ in a finite index subgroup of $\mathrm{U}(\Z)$, then  $\Gamma(f,g)$ has finite index in $\mathrm{O}_\mathrm{Q}(\Z)$. Note that, the criterion of Venkataramana \cite{Ve94} (for $K$ - rank one) and Raghunathan \cite{Rag91} (for $K$ - rank two) are for all absolutely almost simple linear algebraic groups defined over a number field $K$, and we have stated in a way to use it to prove our theorem.

To prove Theorem \ref{maintheorem} for the cases \ref{arithmetic-BS16}, \ref{arithmetic-BS19} and \ref{arithmetic-BS22} of Table \ref{table:arithmetic-Bajpai-Singh} (note that, for all these cases, the corresponding orthogonal groups have $\Q$ - rank two), we use the following criterion of Venkataramana \cite[Theorem 3.5]{Ve87}: If $\Gamma(f,g)\subset\mathrm{O}_\mathrm{Q}(\Z)$ is a Zariski dense subgroup and intersecting the highest and second highest root groups of $\mathrm{SO}_\mathrm{Q}$ non-trivially, then  $\Gamma(f,g)$ has finite index in $\mathrm{O}_\mathrm{Q}(\Z)$.

For an easy reference to the root system, and the structures of the corresponding unipotent subgroups of $\mathrm{SO}_\mathrm{Q}$ (of $\Q$ - rank two), we refer the reader to Remark \ref{methodoftheproof2} (cf. \cite[Section 2.2]{Si15}).

\subsection*{Notation}
Let $f,g$ be a pair of integer coefficient monic polynomials of degree $5$, which have roots of unity as their roots, do not have any common root, form a primitive pair, and satisfy the condition: $f(0)=-1$ and  $g(0)=1$. Let $A, B$ be the companion matrices of $f,g$ respectively, and $\Gamma(f,g)$ be the group generated by $A$ and $B$. Let $C=A^{-1}B$, $I$ be the $5\times 5$ identity matrix; $e_1,e_2,\ldots,e_5$ be the standard basis vectors of $\Q^5$ over $\Q$; and $v$ be the last column vector of $C-I$, that is, $v=(C-I)e_5$. Let $\mathrm{Q}$ be the non-degenerate quadratic form on $\Q^5$, preserved by $\Gamma(f,g)$. Note that, $\mathrm{Q}$ is unique only up to scalars, and its existence follows from Beukers and Heckman \cite{BH89}.

It is clear that $C(v)=-v$ (since $C^2=I$), and hence $v$ is $\mathrm{Q}$ - orthogonal to the vectors $e_1,e_2,e_3,e_4$ (since $\mathrm{Q}(e_i, v)=\mathrm{Q}(C(e_i), C(v))=\mathrm{Q}(e_i, -v)=-\mathrm{Q}(e_i, v),$ for $i=1,2,3,4$); and $\mathrm{Q}(v,e_5)\neq0$ (since $\mathrm{Q}$ is non-degenerate). We may now assume that $\mathrm{Q}(v,e_5)=1$.

It can be checked easily that the set $\{v, Av, A^2v, A^3v, A^4v\}$ (similarly $\{v, Bv, B^2v, B^3v, B^4v\}$) forms a basis of $\Q^5$ (cf. \cite[Lemma 2.1]{Si15}). Therefore, to determine the quadratic form $\mathrm{Q}$ on $\Q^5$, it is enough to compute $\mathrm{Q}(v, A^jv)$, for $j=0,1,2,3,4$. Also, since $v$ is $\mathrm{Q}$ - orthogonal to the vectors $e_1, e_2, e_3, e_4$, and $\mathrm{Q}(v,e_5)=1$ (say), we get 

\begin{equation}\label{computaionofQ}
 \mathrm{Q}(v,A^jv)=\mbox{coefficient of }e_5\mbox{ in }A^jv.
\end{equation}

We compute the signature of the quadratic form $\mathrm{Q}$ by using \cite[Theorem 4.5]{BH89}, which says the following: By renumbering, we may assume that $0\leq\alpha_1\leq\alpha_2\leq\ldots\leq\alpha_5<1$ and $0\leq\beta_1\leq\beta_2\leq\ldots\leq\beta_5<1$. Let $m_j=\#\{k : \beta_k<\alpha_j\}$ for $j=1,2,3,4,5$. Then, the signature $(p,q)$ of the quadratic form $\mathrm{Q}$ is determined by the equation

\begin{equation}\label{computaionofsignature}
\left| p-q\right|=\left|\sum_{j=1}^5(-1)^{j+m_j}\right|. 
\end{equation}

Now using the above equality with the non-degeneracy of the quadratic form $\mathrm{Q}$ (that is, $p+q=5$), one can check easily that all the quadratic forms $\mathrm{Q}$ associated to the pairs of polynomials in Tables \ref{table:arithmetic-Bajpai-Singh} and \ref{table:unknown-Bajpai-Singh} have signature $(3,2)$, that is, the associated orthogonal groups $\mathrm{O}_\mathrm{Q}$ have real rank two. Therefore, the quadratic forms $\mathrm{Q}$ in the cases of Tables \ref{table:arithmetic-Bajpai-Singh} and \ref{table:unknown-Bajpai-Singh}, have two linearly independent isotropic vectors in $\Q^5$ (by the Hasse-Minkowski theorem), which are not $\mathrm{Q}$ - orthogonal, that is, the associated orthogonal groups $\mathrm{O}_\mathrm{Q}$ have $\Q$ - rank $\geq1$.

Therefore, we are able to use the criterion of Raghunathan \cite{Rag91} and Venkataramana \cite{Ve94} (cf. \cite[Theorem 6]{Ve14}) to show the arithmeticity of the hypergeometric groups associated to the pairs of polynomials of Table \ref{table:arithmetic-Bajpai-Singh}, and could not succeed to do the same for the hypergeometric groups associated to the pairs of Table \ref{table:unknown-Bajpai-Singh}.

In the proof of arithmeticity of the hypergeometric groups associated to the pairs of polynomials of Table \ref{table:arithmetic-Bajpai-Singh}, we denote by $X$ the change of basis matrix (cf. second paragraph of Section \ref{proof}), that is, 

\[\e_1=X(e_1),\quad\e_2=X(e_2),\quad\e_3=X(e_3),\quad\e_2^*=X(e_4),\quad\e_1^*=X(e_5).\]\\
We also denote by $a, b$ the matrices $A, B$ respectively, with respect to the new basis, that is, 

\[a=X^{-1}AX,\quad b=X^{-1}BX.\]\\
We compute the quadratic forms $\mathrm{Q}$ on $\Q^5$ by using Equation (\ref{computaionofQ}), and denote by $\mathrm{Q}$ only, the matrix forms of the quadratic forms $\mathrm{Q}$, with respect to the standard basis $\{e_1,e_2,e_3,e_4,e_5\}$ of $\Q^5$, so that the conditions (note that $A^t$, $B^t$ denote the transpose of the respective matrices)

\[A^t\mathrm{Q}A=\mathrm{Q},\quad B^t\mathrm{Q}B=\mathrm{Q}\]
are satisfied.

Let $\mathrm{P}$ be the parabolic $\Q$ - subgroup of $\mathrm{SO}_\mathrm{Q}$ (defined in the third paragraph of Section \ref{proof}), and let $\mathrm{U}$ be the unipotent radical of $\mathrm{P}$. Since $\mathrm{U}(\Z)$ is a free abelian group isomorphic to $\Z^3$, to show that $\Gamma(f,g)\cap\mathrm{U}(\Z)$ has finite index in $\mathrm{U}(\Z)$, it is enough to show that it contains {\it three} linearly independent vectors in $\Z^3$. In the proof below, we denote the three corresponding unipotent elements by $q_{_1}$, $q_{_2}$, and $q_{_3}$ (note that these are the words in $a$ and $b$ which are respectively some conjugates of $A$ and $B$ (see the last paragraph), and hence (with respect to some basis of $\Q^5$) $q_{_1},q_{_2},q_{_3}\in\Gamma(f,g)$); and in some of the $\Q$ - rank two cases, we are able to show only two unipotent elements $q_{_1}$ and $q_{_2}$ inside $\Gamma(f,g)$, which correspond, respectively, to the highest and second highest roots (Remarks \ref{methodoftheproof2}, \ref{methodoftheproof3}, \ref{methodoftheproof4}), and the arithmeticity of $\Gamma(f,g)$ follows; thanks to the criterion \cite[Theorem 3.5]{Ve87} of Venkataramana.

We now note the following remarks:

\begin{rmk}\label{methodoftheproof1}\normalfont
 As explained above, to show the arithmeticity of the hypergeometric groups we need to produce enough unipotent elements in $\Gamma(f,g)$. Note that in the case of symplectic hypergeometric groups (when $n$ is even and $f(0)=g(0)=1$) the matrix $C=A^{-1}B$ itself is a non-trivial unipotent element, and to produce more unipotent elements in $\Gamma(f,g)$, one can take the conjugates of $C$ and then the commutators of those conjugates, by keeping in mind the structure of the required unipotent elements. But in case of the orthogonal hypergeometric groups, the element $C=A^{-1}B$ is not unipotent (since one of the eigen value of $C$ is $-1$ and all other eigen values are $1$) and if we take $C^2$ (to get rid of the eigen value $-1$), it becomes the identity element. Therefore, to find a non-trivial unipotent element in $\Gamma(f,g)$, we take conjugates of $C$ by the elements of $\Gamma(f,g)$, and then try to see the commutators of those conjugates. The idea of finding the first non-trivial unipotent element is completely experiment based and we use Maple extensively to find it in some of the cases. Once we find a non-trivial unipotent element, we take the conjugates of it and then compute the commutators of such conjugates.
 \end{rmk}
 
 \begin{rmk}\label{methodoftheproof2}\normalfont
 For the $\mathbb{Q}$ - rank two cases, the explicit structures of the unipotent subgroups corresponding to the roots of the orthogonal group have been given in \cite[Seciton 2.2]{Si15}. For pedagogical reasons, we summarize it here.
 
 Note that if we consider the maximal torus
 \[\mathrm{T}={\tiny\left\{\begin{pmatrix}
t_1 &0 &0 &0 &0\\
0 &t_2 &0 &0 &0\\
0 &0 &1 &0 &0\\
0 &0 &0 &t_2^{-1} &0\\
0 &0 &0  &0 &t_1^{-1}
\end{pmatrix} : t_i\in\Q^*,\quad \forall\ 1\leq i\leq 2\right\}}\] inside the orthogonal group, and if $\mathbf{t}_i$ is the character of $\mathrm{T}$ defined by
{\tiny\[\begin{pmatrix}
t_1 &0 &0 &0 &0\\
0 &t_2 &0 &0 &0\\
0 &0 &1 &0 &0\\
0 &0 &0 &t_2^{-1} &0\\
0 &0 &0  &0 &t_1^{-1}
\end{pmatrix}\longmapsto t_i,\qquad\mbox{for }i=1,2,\]}then the roots of the orthogonal group are given by $\Phi=\{\mathbf{t}_1, \mathbf{t}_2, \mathbf{t}_1\mathbf{t}_2, \mathbf{t}_1\mathbf{t}_2^{-1},\mathbf{t}_1^{-1},\mathbf{t}_2^{-1}, \mathbf{t}_1^{-1}\mathbf{t}_2^{-1}, \mathbf{t}_2\mathbf{t}_1^{-1}\}.$ If we fix a set of simple roots $\Delta=\{ \mathbf{t}_2, \mathbf{t}_1\mathbf{t}_2^{-1}\},$ then the set of positive roots $\Phi^+=\{\mathbf{t}_2, \mathbf{t}_1\mathbf{t}_2^{-1}, \mathbf{t}_1, \mathbf{t}_1\mathbf{t}_2\}$, the set of negative roots $\Phi^-=\{\mathbf{t}_2^{-1}, \mathbf{t}_2\mathbf{t}_1^{-1}, \mathbf{t}_1^{-1}, \mathbf{t}_1^{-1}\mathbf{t}_2^{-1}\}$, and $\mathbf{t}_1\mathbf{t}_2$, $\mathbf{t}_1$ are respectively the {\it highest}  and  {\it second highest} roots in  $\Phi^+$. 

The unipotent subgroups corresponding to the highest and the second highest roots are respectively given by
\[\mathrm{U}_{\mathbf{t}_1\mathbf{t}_2}={\tiny \left\{\begin{pmatrix}
1 &0 &0 &x &0\\
0 &1 &0 &0 &-\frac{\lambda_1}{\lambda_2}x\\
0 &0 &1 &0 &0\\
0 &0 &0 &1 &0\\
0 &0 &0  &0 &1
\end{pmatrix}: x\in\Q\right\}}\quad\mbox{ and }\quad\mathrm{U}_{\mathbf{t}_1}={\tiny\left\{\begin{pmatrix}
1 &0 &x &0 &-\frac{\lambda_1}{2\lambda_3}x^2\\
0 &1 &0 &0 &0\\
0 &0 &1 &0 &-\frac{\lambda_1}{\lambda_3}x\\
0 &0 &0 &1 &0\\
0 &0 &0  &0 &1
\end{pmatrix}: x\in\Q\right\}}\]
where $\lambda_1=\mathrm{Q}(\epsilon_1,\epsilon_1^*)$, $\lambda_2=\mathrm{Q}(\epsilon_2,\epsilon_2^*)$ and $\lambda_3=\mathrm{Q}(\epsilon_3,\epsilon_3)$, and $\{\epsilon_1,\epsilon_2,\epsilon_3,\epsilon_2^*,\epsilon_1^*\}$ is a basis of $\Q^5$, with respect to which, the matrix associated to the quadratic form $\mathrm{Q}$ has the anti-diagonal form.
\end{rmk}

\begin{rmk}\label{methodoftheproof3}\normalfont
Since $\Gamma(f,g)$ is Zariski dense in the orthogonal group $\mathrm{O}_\mathrm{Q}$ (by \cite[Theorem 6.5]{BH89}), to apply the criterion \cite[Theorem 3.5]{Ve87} in the cases of hypergeometric groups of degree five, for which, the corresponding orthogonal groups $\mathrm{O}_\mathrm{Q}$ have $\Q$ - rank two, it is enough to show that 
\[[\mathrm{U}_{\mathbf{t}_1\mathbf{t}_2}(\Z) : \Gamma(f,g)\cap\mathrm{U}_{\mathbf{t}_1\mathbf{t}_2}(\Z)]\quad\mbox{ and }\quad [\mathrm{U}_{\mathbf{t}_1}(\Z) : \Gamma(f,g)\cap\mathrm{U}_{\mathbf{t}_1}(\Z)]\]are finite. Since  \[[\mathrm{U}(\Z) : \Gamma(f,g)\cap\mathrm{U}(\Z)]<\infty\Rightarrow  [\mathrm{U}_{\mathbf{t}_1\mathbf{t}_2}(\Z) : \Gamma(f,g)\cap\mathrm{U}_{\mathbf{t}_1\mathbf{t}_2}(\Z)]<\infty\]and \[[\mathrm{U}_{\mathbf{t}_1}(\Z) : \Gamma(f,g)\cap\mathrm{U}_{\mathbf{t}_1}(\Z)]<\infty,\]the arithmeticity of the hypergeometric groups in Table \ref{table:arithmetic-Bajpai-Singh}, for which, the associated orthogonal groups have $\mathbb{Q}$ - rank two, also follows from \cite[Theorem 3.5]{Ve87}.
\end{rmk}

\begin{rmk}\label{methodoftheproof4}\normalfont
In the three cases \ref{arithmetic-BS16}, \ref{arithmetic-BS19} and \ref{arithmetic-BS22}, we find the unipotent elements (in $\Gamma(f,g)$) corresponding to the highest and the second highest roots, and finding the third unipotent element in these cases does not seem to be direct compare to the other $20$ cases of Table \ref{table:arithmetic-Bajpai-Singh}, so we do not put much effort in computing the third unipotent element (to show that $[\mathrm{U}(\Z) : \Gamma(f,g)\cap\mathrm{U}(\Z)]<\infty$). Also, once we find that the group $\Gamma(f,g)$ has finite index in the integral orthogonal group (using \cite[Theorem 3.5]{Ve87}) it follows trivially that $\Gamma(f,g)$ contains the unipotent elements corresponding to each roots, and hence  $[\mathrm{U}(\Z) : \Gamma(f,g)\cap\mathrm{U}(\Z)]$ is also finite.
\end{rmk}

\begin{rmk}\label{methodoftheproof5}\normalfont
In some of the remaining $20$ cases (cf. Table \ref{table:unknown-Bajpai-Singh}), for which, the associated quadratic forms $\mathrm{Q}$ have signature $(3,2)$, we find some non-trivial unipotent elements (but not enough to show the arithmeticity), and for the other cases, it is even not possible (for us) to find a non-trivial unipotent element to start with. Thus, we are not able to show the arithmeticity of these groups, and believe that one may apply the ping-pong argument (similar to Brav and Thomas \cite{BT14}) to show the thinness of these groups.
\end{rmk}

We now give a detailed explanation of the computation for finding out the unipotent elements inside the hypergeometric group $\Gamma(f,g)$ in the first case (Case \ref{arithmetic-BS1}) of Table \ref{table:arithmetic-Bajpai-Singh}. Following the similar computation, we find out the unipotent elements for the other cases of Table \ref{table:arithmetic-Bajpai-Singh}.
\subsection{Arithmeticity of Case \ref{arithmetic-BS1}}\label{subsection-BS1}
In this case \[\alpha=\left(0, 0, 0, \frac{1}{3}, \frac{2}{3}\right);\qquad 
\beta=\left(\frac{1}{2},\frac{1}{4},\frac{1}{4},\frac{3}{4},\frac{3}{4}\right)\]
\[f(x)=x^5-2x^4+x^3-x^2+2x-1;\quad g(x)=x^5+x^4+2x^3+2x^2+x+1.\]

Let $A$ and $B$  be the companion  matrices of $f(X)$ and  $g(X)$ respectively, and let $C=A^{-1}B$. Then
\[A={\tiny \begin{pmatrix} \begin {array}{rrrrr} 0&0&0&0&1\\ \noalign{\medskip}1&0&0&0&-2
\\ \noalign{\medskip}0&1&0&0&1\\ \noalign{\medskip}0&0&1&0&-1
\\ \noalign{\medskip}0&0&0&1&2\end {array}
  \end{pmatrix}},\quad  B={\tiny\begin{pmatrix}\begin {array}{rrrrr} 0&0&0&0&-1\\ \noalign{\medskip}1&0&0&0&-
1\\ \noalign{\medskip}0&1&0&0&-2\\ \noalign{\medskip}0&0&1&0&-2
\\ \noalign{\medskip}0&0&0&1&-1\end {array}
 \end{pmatrix}},\quad C={\tiny\begin{pmatrix} \begin {array}{rrrrr} 1&0&0&0&-3\\ \noalign{\medskip}0&1&0&0&-
1\\ \noalign{\medskip}0&0&1&0&-3\\ \noalign{\medskip}0&0&0&1&1
\\ \noalign{\medskip}0&0&0&0&-1
\end {array}
 \end{pmatrix}}.\]

Let $\{e_1, e_2, e_3, e_4, e_5\}$ be the standard basis of $\Q^5$ over $\Q$, and let $v=(C-I)(e_5)=-3e_1-e_2-3e_3+e_4-2e_5$. Then $Bv=2e_1-e_2+3e_3+e_4+3e_5$, $B^2v=-3e_1-e_2-7e_3-3e_4-2e_5$, $B^3v=2e_1-e_2+3e_3-3e_4-e_5$, $B^4v=e_1+3e_2+e_3+5e_4-2e_5$; and hence by the Equation (\ref{computaionofQ}) it follows that, with respect to the basis $\{v, Bv, B^2v, B^3v, B^4v\}$, the matrix of the quadratic form (up to scalar multiplication) preserved by the hypergeometric group $\Gamma(f,g)$ is given by
\[\mathrm{Q}={\tiny\begin{pmatrix}\begin {array}{rrrrr} -2&3&-2&-1&-2\\ \noalign{\medskip}3&-2&3
&-2&-1\\ \noalign{\medskip}-2&3&-2&3&-2\\ \noalign{\medskip}-1&-2&3&-2
&3\\ \noalign{\medskip}-2&-1&-2&3&-2
\end {array}
\end{pmatrix}}.\]

We now let $Z$ be the change of basis matrix from the basis $\{v, Bv, B^2v, B^3v, B^4v\}$ to the standard basis $\{e_1, e_2, e_3, e_4, e_5\}$. Then, a computation shows that
\[Z={\tiny\begin{pmatrix}
     \begin {array}{rrrrr} -{{5}/{32}}&-{{11}/{32}}&{{3}/{32}}&{5}/{32}&-{{5}/{32}}\\ \noalign{\medskip}{1}/{2}&-{1}/{2}&-{1}/{4}&{1}/{4}&0\\ \noalign{\medskip}{{7}/{16}}&-{3}/{16}&-{5}/{16}&{1}/{16}&-{1}/{16}\\ \noalign{\medskip}{1}/{4}&-{1}/{4}&0&0&-{1}/{4}
\\ \noalign{\medskip}{{11}/{32}}&-{{3}/{32}}&-{{5}/{32}
}&{{5}/{32}}&-{{5}/{32}}\end {array}
    \end{pmatrix}}
\]and, with respect to the standard basis $\{e_1, e_2, e_3, e_4, e_5\}$, the matrix (up to scalar multiplication) of the quadratic form preserved by the hypergeometric group $\Gamma(f,g)$ is given by
\[-32Z^tQZ={\tiny\begin{pmatrix}
            \begin {array}{rrrrr} 5&-5&-3&11&5\\ \noalign{\medskip}-5&5&-5
&-3&11\\ \noalign{\medskip}-3&-5&5&-5&-3\\ \noalign{\medskip}11&-3&-5&
5&-5\\ \noalign{\medskip}5&11&-3&-5&5\end {array}
           \end{pmatrix}}\](here we multiply by $-32$, just to get rid of the denominators) where $Z^t$ denotes the transpose of the matrix $Z$. Thus, we find that \[A^t(-32Z^tQZ)A=(-32Z^tQZ)=B^t(-32Z^tQZ)B.\]
           
After computing the isotropic vectors of the quadratic form $Q$, we change the standard basis $\{e_1,e_2,e_3, e_4,e_5\}$ to the basis $\{\e_1,\e_2,\e_3,\e_2^*,\e_1^*\}$, for which, the change of basis matrix is given by
\[X={\tiny\begin{pmatrix}\begin {array}{rrrrr} 4&0&0&0&-1\\ \noalign{\medskip}4&0&-2&6&
2\\ \noalign{\medskip}4&4&-2&2&-2\\ \noalign{\medskip}4&4&-6&10&0
\\ \noalign{\medskip}0&0&2&-2&-3\end {array}\end{pmatrix}}.\]With respect to the new basis, the matrices of $\mathrm{Q}$, $A$ and $B$ are $X^t\mathrm{Q}X$, $a=X^{-1}AX$ and $b=X^{-1}BX$; which are respectively
{\tiny\[\begin{pmatrix}\begin {array}{rrrrr} 0&0&0&0&-128\\ \noalign{\medskip}0&0&0&-
128&0\\ \noalign{\medskip}0&0&64&0&0\\ \noalign{\medskip}0&-128&0&0&0
\\ \noalign{\medskip}-128&0&0&0&0\end {array}\end{pmatrix},\quad \begin{pmatrix}\begin {array}{rrrrr} 0&-1/2&1/2&-1/2&0\\ \noalign{\medskip}2&
-1&-1&3&2\\ \noalign{\medskip}4&0&-3&6&2\\ \noalign{\medskip}2&1&-2&3&
1/2\\ \noalign{\medskip}0&-2&0&0&3\end {array}\end{pmatrix},\quad \begin{pmatrix}\begin {array}{rrrrr} 0&-1/2&0&0&3/4\\ \noalign{\medskip}2&-1&
-1&3&2\\ \noalign{\medskip}4&0&-3&6&2\\ \noalign{\medskip}2&1&-2&3&1/2
\\ \noalign{\medskip}0&-2&2&-2&0\end{array}\end{pmatrix}.\]}

It is clear from the above computation that the $\Q$ - rank of the orthogonal group $\mathrm{SO}_\mathrm{Q}$ is two.
If we denote by
\[w_{_1} = ab^{-1}ab,\quad w_{_2} = a^{-1}b^4a,\quad w_{_3} = w_{_1}^2w_{_2}^{-1},\quad w_{_4} = (w_{_3} w_{_1}^{-1})^{2} w_{_3},\]
then \[w_1={\tiny\begin{pmatrix}
                  \begin {array}{rrrrr} -1&-1&2&-3&-1\\ \noalign{\medskip}0&-1&2
&-4&-1\\ \noalign{\medskip}0&0&1&-4&0\\ \noalign{\medskip}0&0&0&-1&1
\\ \noalign{\medskip}0&0&0&0&-1\end {array}
                 \end{pmatrix}},\quad w_2={\tiny\begin{pmatrix}
                  \begin {array}{rrrrr}1&2&0&0&0\\ \noalign{\medskip}0&1&0&0&0
\\ \noalign{\medskip}0&0&1&0&0\\ \noalign{\medskip}0&0&0&1&-2
\\ \noalign{\medskip}0&0&0&0&1
\end {array}
                 \end{pmatrix}},\]
                 
                 \[ w_3={\tiny\begin{pmatrix}
                  \begin {array}{rrrrr} 1&0&-2&2&4\\ \noalign{\medskip}0&1&0&0&-
2\\ \noalign{\medskip}0&0&1&0&-4\\ \noalign{\medskip}0&0&0&1&0
\\ \noalign{\medskip}0&0&0&0&1
\end {array}
                 \end{pmatrix}},\quad w_4={\tiny\begin{pmatrix}
                  \begin {array}{rrrrr} 1&-2&0&-4&-8\\ \noalign{\medskip}0&1&0&0
&4\\ \noalign{\medskip}0&0&1&0&0\\ \noalign{\medskip}0&0&0&1&2
\\ \noalign{\medskip}0&0&0&0&1
\end {array}
                 \end{pmatrix}},\] and the final three unipotent elements are
\[q_{_1}=w_{_4} w_{_2},\quad q_{_2}=w_{_3}^{2} q_{_1},\quad q_{_3}=w_{_2},\] which are respectively
{\tiny\[\begin{pmatrix}\begin {array}{rrrrr} 1&0&0&-4&0\\ \noalign{\medskip}0&1&0&0&4
\\ \noalign{\medskip}0&0&1&0&0\\ \noalign{\medskip}0&0&0&1&0
\\ \noalign{\medskip}0&0&0&0&1\end {array}\end{pmatrix},\quad \begin{pmatrix}\begin {array}{rrrrr} 1&0&-4&0&16\\ \noalign{\medskip}0&1&0&0&0
\\ \noalign{\medskip}0&0&1&0&-8\\ \noalign{\medskip}0&0&0&1&0
\\ \noalign{\medskip}0&0&0&0&1\end {array}\end{pmatrix},\quad \begin{pmatrix}\begin {array}{rrrrr} 1&2&0&0&0\\ \noalign{\medskip}0&1&0&0&0
\\ \noalign{\medskip}0&0&1&0&0\\ \noalign{\medskip}0&0&0&1&-2
\\ \noalign{\medskip}0&0&0&0&1\end {array}\end{pmatrix}.\]}It is clear from the above computation that the group generated by $q_1, q_2, q_3$ is a subgroup of $\Gamma(f,g)\cap\mathrm{U}(\Z)$, and has finite index in $\mathrm{U}(\Z)$, the integer points of the unipotent radical $\mathrm{U}$ of the parabolic subgroup $\mathrm{P}$, which preserves the flag \[\{0\}\subset\Q\e_1\subset\Q\e_1\oplus\Q\e_2\oplus\Q\e_3\oplus\Q\e_2^*\subset\Q\e_1\oplus\Q\e_2\oplus\Q\e_3\oplus\Q\e_2^*\oplus\Q\e_1^*.\] Therefore $\Gamma(f,g)\cap\mathrm{U}(\Z)$ is of finite index in $\mathrm{U}(\Z)$. Since $\Gamma(f,g)$ is also Zariski dense in the orthogonal group $\mathrm{O}_\mathrm{Q}$ (by Beukers and Heckman \cite{BH89}), the arithmeticity of $\Gamma(f,g)$ follows from Venkataramana \cite{Ve94} (cf. \cite{Rag91}, \cite[Theorem 6]{Ve14}).

\begin{rmk}\normalfont
 We follow the same notation as that of Section \ref{subsection-BS1} for the remaining cases, and the arithmeticity of $\Gamma(f,g)$ follows by the similar arguments used in Section \ref{subsection-BS1}. Also, in the proof of the arithmeticity of the remaining cases, we denote the commutator of two elements $s,t$ (say) in $\Gamma(f,g)$ by $[s,t]$ which is $sts^{-1}t^{-1}$.
\end{rmk}

\subsection{Arithmeticity of Case \ref{arithmetic-BS2}}\label{subsection-BS2}
In this case \[\alpha=\left(0, 0, 0, \frac{1}{3}, \frac{2}{3}\right);\qquad 
\beta=\left(\frac{1}{2},\frac{1}{4},\frac{3}{4},\frac{1}{6},\frac{5}{6}\right)\]
\[f(x)=x^5-2x^4+x^3-x^2+2x-1;\quad g(x)=x^5 +x^3+x^2+1.\]

The matrix (up to scalar multiplication) of the quadratic form $\mathrm{Q}$, 
with respect to the standard basis $\{e_1, e_2, e_3, e_4, e_5\}$, and  the change of basis matrix $X$, are
{\tiny\[ \mathrm{Q}=\begin{pmatrix}\begin {array}{rrrrr}1&3&1&-5&-7\\ \noalign{\medskip}3&1&3&1&
-5\\ \noalign{\medskip}1&3&1&3&1\\ \noalign{\medskip}-5&1&3&1&3
\\ \noalign{\medskip}-7&-5&1&3&1
\end {array}
\end{pmatrix};\quad X=\begin{pmatrix}\begin {array}{rrrrr} 0&2&-2&-1/2&1\\ \noalign{\medskip}1&0&0&0
&-1/2\\ \noalign{\medskip}0&4&-2&0&2\\ \noalign{\medskip}-1&0&-2&-3&-1
/2\\ \noalign{\medskip}0&2&-2&-1/2&0
\end {array}\end{pmatrix},\]}and the matrices of $\mathrm{Q}$, $A$ and $B$, with respect to the new basis, are $X^t\mathrm{Q}X$, $a=X^{-1}AX$ and $b=X^{-1}BX$; which are respectively
{\tiny\[\begin{pmatrix}\begin {array}{rrrrr} 0&0&0&0&8\\ \noalign{\medskip}0&0&0&-16&0
\\ \noalign{\medskip}0&0&-16&0&0\\ \noalign{\medskip}0&-16&0&0&0
\\ \noalign{\medskip}8&0&0&0&0
\end {array}\end{pmatrix},\quad \begin{pmatrix}\begin {array}{rrrrr} 1/2&-3&4&9/4&5/4\\ \noalign{\medskip}1/4
&1/2&-1&-{9}/{8}&-1/8\\ \noalign{\medskip}1&-2&3&3/2&1/2
\\ \noalign{\medskip}-1&2&-4&-5/2&-3/2\\ \noalign{\medskip}1&-2&4&7/2&
1/2
\end {array}\end{pmatrix},\quad \begin{pmatrix}\begin {array}{rrrrr} 1/2&1&0&5/4&5/4\\ \noalign{\medskip}1/4&
3/2&-2&-{{11}/{8}}&-1/8\\ \noalign{\medskip}1&2&-1&1/2&1/2
\\ \noalign{\medskip}-1&-2&0&-3/2&-3/2\\ \noalign{\medskip}1&-2&4&7/2&
1/2
\end{array}\end{pmatrix}.\]}

It is clear from the above computation that the $\Q$ - rank of the orthogonal group $\mathrm{SO}_\mathrm{Q}$ is two. If we denote by 
\[w_{_1}=[a^{-1}\,, b],\quad w_{_2}=(aba)^2 (bab)^{-1},\quad w_{_3}=[w_{_2}^{2}\,, w_{_1}],\]then
\[w_1={\tiny\begin{pmatrix}
                  \begin {array}{rrrrr} 1&4&-4&-1&2\\ \noalign{\medskip}0&1&0&0&
-1/2\\ \noalign{\medskip}0&0&1&0&-2\\ \noalign{\medskip}0&0&0&1&2
\\ \noalign{\medskip}0&0&0&0&1
\end {array}
                 \end{pmatrix}},\quad w_2={\tiny\begin{pmatrix}
                  \begin {array}{rrrrr}-1&0&0&1&0\\ \noalign{\medskip}0&-1&-1&1
/2&1/2\\ \noalign{\medskip}-2&0&1&1&0\\ \noalign{\medskip}0&0&0&-1&0
\\ \noalign{\medskip}-4&0&4&0&-1
\end {array}
                 \end{pmatrix}},\quad w_3={\tiny\begin{pmatrix}
                  \begin {array}{rrrrr} 1&0&-8&-16&16\\ \noalign{\medskip}0&1&0&0
&-8\\ \noalign{\medskip}0&0&1&0&-4\\ \noalign{\medskip}0&0&0&1&0
\\ \noalign{\medskip}0&0&0&0&1
\end {array}
                 \end{pmatrix}},\]and the final three unipotent elements are
\[q_{_1}=[ w_{_2}^{2} \,, w_{_3}],\quad q_{_2}=w_{_3} q_{_1}^{-1},\quad q_{_3}=(w_{_1}^{-4} q_{_2}^{2})^{4} q_{_1},\]which are respectively
{\tiny\[\begin{pmatrix}\begin {array}{rrrrr} 1&0&0&-16&0\\ \noalign{\medskip}0&1&0&0&
-8\\ \noalign{\medskip}0&0&1&0&0\\ \noalign{\medskip}0&0&0&1&0
\\ \noalign{\medskip}0&0&0&0&1
\end {array}\end{pmatrix},\quad \begin{pmatrix}\begin {array}{rrrrr} 1&0&-8&0&16\\ \noalign{\medskip}0&1&0&0&0
\\ \noalign{\medskip}0&0&1&0&-4\\ \noalign{\medskip}0&0&0&1&0
\\ \noalign{\medskip}0&0&0&0&1
\end {array}\end{pmatrix},\quad \begin{pmatrix}\begin {array}{rrrrr}1&-64&0&0&0\\ \noalign{\medskip}0&1&0&0&0
\\ \noalign{\medskip}0&0&1&0&0\\ \noalign{\medskip}0&0&0&1&-32
\\ \noalign{\medskip}0&0&0&0&1
\end {array}\end{pmatrix}.\]}

\subsection{Arithmeticity of Case \ref{arithmetic-BS3}}\label{subsection-BS3}
In this case \[\alpha=\left(0, 0, 0, \frac{1}{3}, \frac{2}{3}\right);\qquad 
\beta=\left(\frac{1}{2},\frac{1}{6},\frac{1}{6},\frac{5}{6},\frac{5}{6}\right)\]
\[f(x)=x^5-2x^4+x^3-x^2+2x-1;\quad g(x)=x^5-x^4+x^3+x^2-x+1.\]

The matrix (up to scalar multiplication) of the quadratic form $\mathrm{Q}$, 
with respect to the standard basis $\{e_1, e_2, e_3, e_4, e_5\}$, and  the change of basis matrix $X$, are
{\tiny\[ \mathrm{Q}=\begin{pmatrix}\begin {array}{rrrrr}3&3&1&-3&-7\\ \noalign{\medskip}3&3&3&1&
-3\\ \noalign{\medskip}1&3&3&3&1\\ \noalign{\medskip}-3&1&3&3&3
\\ \noalign{\medskip}-7&-3&1&3&3
\end {array}
\end{pmatrix};\quad X=\begin{pmatrix}\begin {array}{rrrrr} -2&0&0&0&1\\ \noalign{\medskip}0&-2&8&-1
&-3\\ \noalign{\medskip}0&0&-8&2&2\\ \noalign{\medskip}-2&0&8&-4&-1
\\ \noalign{\medskip}0&-2&0&1&1
\end {array}\end{pmatrix},\]}and the matrices of $\mathrm{Q}$, $A$ and $B$, with respect to the new basis, are $X^t\mathrm{Q}X$, $a=X^{-1}AX$ and $b=X^{-1}BX$; which are respectively
{\tiny\[\begin{pmatrix}\begin {array}{rrrrr} 0&0&0&0&16\\ \noalign{\medskip}0&0&0&16&0
\\ \noalign{\medskip}0&0&-64&0&0\\ \noalign{\medskip}0&16&0&0&0
\\ \noalign{\medskip}16&0&0&0&0
\end {array}\end{pmatrix},\quad \begin{pmatrix}\begin {array}{rrrrr} 0&0&0&-1/2&1/2\\ \noalign{\medskip}1&1&-
4&1&1\\ \noalign{\medskip}0&0&-1&0&1\\ \noalign{\medskip}0&0&0&0&1
\\ \noalign{\medskip}0&-2&0&0&2
\end {array}\end{pmatrix},\quad \begin{pmatrix}\begin {array}{rrrrr} 0&-1&0&0&1\\ \noalign{\medskip}1&1&-4&1&
1\\ \noalign{\medskip}0&0&-1&0&1\\ \noalign{\medskip}0&0&0&0&1
\\ \noalign{\medskip}0&0&0&-1&1
\end{array}\end{pmatrix}.\]}

It is clear from the above computation that the $\Q$ - rank of the orthogonal group $\mathrm{SO}_\mathrm{Q}$ is two. If we denote by \[d=b a^{-1},\quad w_{_1}=a^{-1}b^{6} a,\quad w_{_2}=a b^{-6} a^{-1},\quad w_{_3}=(a b^{-4})^{2},\quad w_{_4}=(a^{-1} b^{4})^{2},\]
 \[w_{_5}=[a^{-1}, b] a^{-2},\quad w_{_6}=d w_{_1} d,\quad w_{_7}=w_{_6}^{-4} w_{_5} w_{_1} w_{_5}^{-1} w_{_2},\]then
\[d={\tiny\begin{pmatrix}
                  \begin {array}{rrrrr}0&0&0&0&1/2\\ \noalign{\medskip}0&1&0&0&0
\\ \noalign{\medskip}0&0&1&0&0\\ \noalign{\medskip}0&0&0&1&0
\\ \noalign{\medskip}2&0&0&0&0
\end {array}
                 \end{pmatrix}},\quad w_1={\tiny\begin{pmatrix}
                  \begin {array}{rrrrr}1&4&0&0&0\\ \noalign{\medskip}0&1&0&0&0
\\ \noalign{\medskip}0&0&1&0&0\\ \noalign{\medskip}0&0&0&1&-4
\\ \noalign{\medskip}0&0&0&0&1
\end {array}
                 \end{pmatrix}},\quad w_2={\tiny\begin{pmatrix}
                  \begin {array}{rrrrr}1&0&0&0&0\\ \noalign{\medskip}4&1&0&0&0
\\ \noalign{\medskip}0&0&1&0&0\\ \noalign{\medskip}0&0&0&1&0
\\ \noalign{\medskip}0&0&0&-4&1
\end {array}
                 \end{pmatrix}},\] \[w_3={\tiny\begin{pmatrix}
                  \begin {array}{rrrrr}1&0&0&-1&0\\ \noalign{\medskip}2&1&0&-2&
1\\ \noalign{\medskip}0&0&1&0&0\\ \noalign{\medskip}0&0&0&1&0
\\ \noalign{\medskip}0&0&0&-2&1
\end {array}
                 \end{pmatrix}},\quad w_4={\tiny\begin{pmatrix}
                  \begin {array}{rrrrr}1&2&0&1&-2\\ \noalign{\medskip}0&1&0&0&-
1\\ \noalign{\medskip}0&0&1&0&0\\ \noalign{\medskip}0&0&0&1&-2
\\ \noalign{\medskip}0&0&0&0&1
\end {array}
                 \end{pmatrix}},\quad w_5={\tiny\begin{pmatrix}
                  \begin {array}{rrrrr}0&1&0&0&0\\ \noalign{\medskip}-1&1&0&0&0
\\ \noalign{\medskip}-2&0&1&0&0\\ \noalign{\medskip}-8&-4&8&0&-1
\\ \noalign{\medskip}4&4&-8&1&1
\end {array}
                 \end{pmatrix}},\] \[w_6={\tiny\begin{pmatrix}
                  \begin {array}{rrrrr} 1&0&0&0&0\\ \noalign{\medskip}0&1&0&0&0
\\ \noalign{\medskip}0&0&1&0&0\\ \noalign{\medskip}-8&0&0&1&0
\\ \noalign{\medskip}0&8&0&0&1
\end {array}
                 \end{pmatrix}},\quad w_7={\tiny\begin{pmatrix}
                  \begin {array}{rrrrr}1&0&0&0&0\\ \noalign{\medskip}0&1&0&0&0
\\ \noalign{\medskip}-8&0&1&0&0\\ \noalign{\medskip}0&0&0&1&0
\\ \noalign{\medskip}128&0&-32&0&1

\end {array}
                 \end{pmatrix}},\]and the final three unipotent elements are
\[q_{_1}=w_{_4}^{2} w_{_1}^{-1},\quad q_{_2}=d w_{_7} d,\quad q_{_3}=w_{_1},\]which are respectively
{\tiny\[\begin{pmatrix}\begin {array}{rrrrr}1&0&0&2&0\\ \noalign{\medskip}0&1&0&0&-2
\\ \noalign{\medskip}0&0&1&0&0\\ \noalign{\medskip}0&0&0&1&0
\\ \noalign{\medskip}0&0&0&0&1
\end {array}\end{pmatrix},\quad \begin{pmatrix}\begin {array}{rrrrr} 1&0&-16&0&32\\ \noalign{\medskip}0&1&0&0
&0\\ \noalign{\medskip}0&0&1&0&-4\\ \noalign{\medskip}0&0&0&1&0
\\ \noalign{\medskip}0&0&0&0&1
\end {array}\end{pmatrix},\quad \begin{pmatrix}\begin {array}{rrrrr}1&4&0&0&0\\ \noalign{\medskip}0&1&0&0&0
\\ \noalign{\medskip}0&0&1&0&0\\ \noalign{\medskip}0&0&0&1&-4
\\ \noalign{\medskip}0&0&0&0&1
\end {array}\end{pmatrix}.\]}

\subsection{Arithmeticity of Case \ref{arithmetic-BS4}}\label{subsection-BS4}
In this case \[\alpha=\left(0, 0, 0, \frac{1}{4}, \frac{3}{4}\right);\qquad 
\beta=\left(\frac{1}{2},\frac{1}{6},\frac{1}{6},\frac{5}{6},\frac{5}{6}\right)\]
\[f(x)=x^5-3x^4+4x^3-4x^2+3x-1;\quad g(x)=x^5-x^4+x^3+x^2-x+1.\]

The matrix (up to scalar multiplication) of the quadratic form $\mathrm{Q}$, 
with respect to the standard basis $\{e_1, e_2, e_3, e_4, e_5\}$, and  the change of basis matrix $X$, are
{\tiny\[\mathrm{Q}=\begin{pmatrix}\begin {array}{rrrrr}7&9&7&-7&-25\\ \noalign{\medskip}9&7&9&7
&-7\\ \noalign{\medskip}7&9&7&9&7\\ \noalign{\medskip}-7&7&9&7&9
\\ \noalign{\medskip}-25&-7&7&9&7
\end {array}
\end{pmatrix};\quad X=\begin{pmatrix}\begin {array}{rrrrr} 0&2&0&1&-2\\ \noalign{\medskip}-1&-2&0&-
1&3\\ \noalign{\medskip}0&4&4&3&-6\\ \noalign{\medskip}1&-6&-4&-2&4
\\ \noalign{\medskip}0&2&0&1&-1
\end {array}\end{pmatrix},\]}and the matrices of $\mathrm{Q}$, $A$ and $B$, with respect to the new basis, are $X^t\mathrm{Q}X$, $a=X^{-1}AX$ and $b=X^{-1}BX$; which are respectively
{\tiny\[\begin{pmatrix}\begin {array}{rrrrr}0&0&0&0&16\\ \noalign{\medskip}0&0&0&-32
&0\\ \noalign{\medskip}0&0&-64&0&0\\ \noalign{\medskip}0&-32&0&0&0
\\ \noalign{\medskip}16&0&0&0&0
\end {array}\end{pmatrix},\quad \begin{pmatrix}\begin {array}{rrrrr} 1&0&-4&1&2\\ \noalign{\medskip}1/2&0&-2&0
&1\\ \noalign{\medskip}0&0&-1&0&1\\ \noalign{\medskip}1&-2&-4&1&1
\\ \noalign{\medskip}1&-2&-4&0&2
\end {array}\end{pmatrix},\quad \begin{pmatrix}\begin {array}{rrrrr} 1&-4&-4&-1&4\\ \noalign{\medskip}1/2&-1&
-2&-1/2&3/2\\ \noalign{\medskip}0&0&-1&0&1\\ \noalign{\medskip}1&-4&-4
&0&2\\ \noalign{\medskip}1&-2&-4&0&2
\end{array}\end{pmatrix}.\]}

It is clear from the above computation that the $\Q$ - rank of the orthogonal group $\mathrm{SO}_\mathrm{Q}$ is two. If we denote by \[c=a^{-1} b,\quad d=[c, a],\quad e=[a^{-1}, b^{2}],\quad s=(c a c)^{4},\quad w_{_1}=e^{-1} d e^{-1} d^{-1},\quad w_{_2}=d [e, c] d^{-1},\] 
\[w_{_3}=w_{_1}^{-2} d^{-1} w_{_2} s w_{_2}^{-1} d,\quad w_{_4}=w_{_3}^{-1} s d^{-4},\quad w_{_5}=[ w_{_4} , d ],\quad w_{_6}=a^{4} w_{_4} d w_{_4}^{-1} d^{7} a^{-4}\]
\[w_{_7}=d^{-82} (w_{_4} d w_{_4}^{-1})^{2},\quad w_{_8}=w_{_6}^{2}(w_{_4} d w_{_4}^{-1} d^{7})^{-10},\]then the final three unipotent elements are
\[q_{_1}=[ w_{_6}, w_{_8}],\quad q_{_2}=w_{_7}^{-240} q_{_3},\quad q_{_3}=(w_{_5}^{-2} w_{_7})^{288} q_{_1},\]
which are respectively
{\tiny\[\begin{pmatrix}\begin {array}{rrrrr}1&0&0&-23040&0\\ \noalign{\medskip}0&1&0
&0&-11520\\ \noalign{\medskip}0&0&1&0&0\\ \noalign{\medskip}0&0&0&1&0
\\ \noalign{\medskip}0&0&0&0&1
\end {array}\end{pmatrix},\quad \begin{pmatrix}\begin {array}{rrrrr}1&0&23040&0&66355200
\\ \noalign{\medskip}0&1&0&0&0\\ \noalign{\medskip}0&0&1&0&5760
\\ \noalign{\medskip}0&0&0&1&0\\ \noalign{\medskip}0&0&0&0&1
\end {array}\end{pmatrix},\quad \begin{pmatrix}\begin {array}{rrrrr}1&46080&0&0&0\\ \noalign{\medskip}0&1&0&0
&0\\ \noalign{\medskip}0&0&1&0&0\\ \noalign{\medskip}0&0&0&1&23040
\\ \noalign{\medskip}0&0&0&0&1

\end {array}\end{pmatrix}.\]}

\subsection{Arithmeticity of Case \ref{arithmetic-BS5}}\label{subsection-BS5}
In this case \[\alpha=\left(0, 0, 0, \frac{1}{6}, \frac{5}{6}\right);\qquad 
\beta=\left(\frac{1}{2},\frac{1}{4},\frac{1}{4},\frac{3}{4},\frac{3}{4}\right)\]
\[f(x)=x^5-4x^4+7x^3-7x^2+4x-1;\quad g(x)=x^5 +x^4+2x^3+2x^2+x+1.\]

The matrix (up to scalar multiplication) of the quadratic form $\mathrm{Q}$, 
with respect to the standard basis $\{e_1, e_2, e_3, e_4, e_5\}$, and  the change of basis matrix $X$, are
{\tiny\[\mathrm{Q}=\begin{pmatrix}\begin {array}{rrrrr}37&11&-35&-37&37\\ \noalign{\medskip}11&
37&11&-35&-37\\ \noalign{\medskip}-35&11&37&11&-35
\\ \noalign{\medskip}-37&-35&11&37&11\\ \noalign{\medskip}37&-37&-35&
11&37
\end {array}
\end{pmatrix};\quad X=\begin{pmatrix}\begin {array}{rrrrr}4&0&-24&-10&11\\ \noalign{\medskip}4&2&-
24&-9&6\\ \noalign{\medskip}4&2&-48&-19&18\\ \noalign{\medskip}4&2&-24
&-7&8\\ \noalign{\medskip}0&2&-24&-11&9
\end {array}\end{pmatrix},\]}and the matrices of $\mathrm{Q}$, $A$ and $B$, with respect to the new basis, are $X^t\mathrm{Q}X$, $a=X^{-1}AX$ and $b=X^{-1}BX$; which are respectively
{\tiny\[\begin{pmatrix}\begin {array}{rrrrr}0&0&0&0&-384\\ \noalign{\medskip}0&0&0&
192&0\\ \noalign{\medskip}0&0&2304&0&0\\ \noalign{\medskip}0&192&0&0&0
\\ \noalign{\medskip}-384&0&0&0&0
\end {array}\end{pmatrix},\quad \begin{pmatrix}\begin {array}{rrrrr} 0&-6&66&{{121}/{4}}&-{{47}/{2}}\\ \noalign{\medskip}2&-10&108&{{105}/{2}}&-39
\\ \noalign{\medskip}0&-2&23&11&-8\\ \noalign{\medskip}0&0&0&0&-1
\\ \noalign{\medskip}0&-2&24&12&-9
\end {array}\end{pmatrix},\quad \begin{pmatrix}\begin {array}{rrrrr}0&-1/2&0&0&5/4\\ \noalign{\medskip}2&0&-
12&-5/2&6\\ \noalign{\medskip}0&0&-1&0&1\\ \noalign{\medskip}0&0&0&0&-
1\\ \noalign{\medskip}0&0&0&1&0
\end{array}\end{pmatrix}.\]}

It is clear from the above computation that the $\Q$ - rank of the orthogonal group $\mathrm{SO}_\mathrm{Q}$ is two. If we denote by \[w_{_1}=b^{4},\quad w_{_2}=a^{-1} b^{4} a,\quad 
w_{_3}=w_{_2}w_{_1}^{-12},\quad w_{_4}=[a, b^{4}],\]
 \[w_{_5}=(ab)^{-1} b^{-2} (ab)^{-1} b^2,\quad 
w_{_6}=w_{_1}^{-1}w_{_2}w_{_3}^{-1} w_{_4}^{-1} w_{_5}^{-1},\]then the final three unipotent elements are
\[q_{_1}=w_{_1},\quad q_{_2}=w_{_5}^{-1} w_{_4}^{-1} w_{_2} w_{_6}^{-1} 
w_{_2}^{-1} q_{_1}^{-277},\quad
q_{_3}=w_{_2}^{-32} q_{_1}^{384} q_{_2}^{8},\]
which are respectively
{\tiny\[\begin{pmatrix}\begin {array}{rrrrr}1&0&0&1&0\\ \noalign{\medskip}0&1&0&0&2
\\ \noalign{\medskip}0&0&1&0&0\\ \noalign{\medskip}0&0&0&1&0
\\ \noalign{\medskip}0&0&0&0&1
\end {array}\end{pmatrix},\quad \begin{pmatrix}\begin {array}{rrrrr}1&0&96&0&768\\ \noalign{\medskip}0&1&0&0
&0\\ \noalign{\medskip}0&0&1&0&16\\ \noalign{\medskip}0&0&0&1&0
\\ \noalign{\medskip}0&0&0&0&1
\end {array}\end{pmatrix},\quad \begin{pmatrix}\begin {array}{rrrrr}1&64&0&0&0\\ \noalign{\medskip}0&1&0&0&0
\\ \noalign{\medskip}0&0&1&0&0\\ \noalign{\medskip}0&0&0&1&128
\\ \noalign{\medskip}0&0&0&0&1
\end {array}\end{pmatrix}.\]}

\subsection{Arithmeticity of Case \ref{arithmetic-BS6}}\label{subsection-BS6}
In this case \[\alpha=\left(0, \frac{1}{3},\frac{1}{3}, \frac{2}{3}, \frac{2}{3}\right);\qquad \beta=\left(\frac{1}{2},\frac{1}{4},\frac{3}{4},\frac{1}{6},\frac{5}{6}\right)\]
\[f(x)=x^5+x^4+x^3-x^2-x-1;\quad g(x)=x^5 +x^3+x^2+1.\]

The matrix (up to scalar multiplication) of the quadratic form $\mathrm{Q}$, 
with respect to the standard basis $\{e_1, e_2, e_3, e_4, e_5\}$, and  the change of basis matrix $X$, are
{\tiny\[\mathrm{Q}=\begin{pmatrix}\begin {array}{rrrrr}1&0&-2&1&2\\ \noalign{\medskip}0&1&0&-2&
1\\ \noalign{\medskip}-2&0&1&0&-2\\ \noalign{\medskip}1&-2&0&1&0
\\ \noalign{\medskip}2&1&-2&0&1
\end {array}
\end{pmatrix};\quad X=\begin{pmatrix}\begin {array}{rrrrr}-2&0&4&0&-1\\ \noalign{\medskip}0&2&0&1&
1\\ \noalign{\medskip}-4&0&4&0&0\\ \noalign{\medskip}-2&4&0&0&1
\\ \noalign{\medskip}0&2&-4&-1&3
\end {array}\end{pmatrix},\]}and the matrices of $\mathrm{Q}$, $A$ and $B$, with respect to the new basis, are $X^t\mathrm{Q}X$, $a=X^{-1}AX$ and $b=X^{-1}BX$; which are respectively
{\tiny\[\begin{pmatrix}\begin {array}{rrrrr}0&0&0&0&8\\ \noalign{\medskip}0&0&0&-8&0
\\ \noalign{\medskip}0&0&-16&0&0\\ \noalign{\medskip}0&-8&0&0&0
\\ \noalign{\medskip}8&0&0&0&0
\end {array}\end{pmatrix},\quad \begin{pmatrix}\begin {array}{rrrrr}0&1&-2&-1&1\\ \noalign{\medskip}-1&-1&2&0
&-1\\ \noalign{\medskip}0&2&-3&-1&2\\ \noalign{\medskip}0&0&0&0&1
\\ \noalign{\medskip}0&4&-4&-1&3
\end {array}\end{pmatrix},\quad \begin{pmatrix}\begin {array}{rrrrr}0&0&0&-1/2&-1/2\\ \noalign{\medskip}-1&-
1&2&0&-1\\ \noalign{\medskip}0&0&1&0&-1\\ \noalign{\medskip}0&0&0&0&1
\\ \noalign{\medskip}0&2&0&0&0
\end{array}\end{pmatrix}.\]}

It is clear from the above computation that the $\Q$ - rank of the orthogonal group $\mathrm{SO}_\mathrm{Q}$ is two. If we denote by \[c=a^{-1} b,\ %d &=& a^{-3},\quad
e=[ a^{2}, b^{2}],\ u=[ a^{-1}, b^{2}],\quad w=[a^{-1}, b^{4}],\quad
%r & = & (a b^{-4} )^{2},\quad
t=w^{-1} u e^{-1},\ w_{_1}=(t^{2} a^{3})^{-2} t^{4} (a^{3} (b^{4} a^{-1})^{3})^{2},\]
%z & = & b^4,\quad
\[w_{_2}=w_{_1} (t^{2} a^{3})^{-2},\quad
w_{_3}=w_{_2} ( w_{_2}^{-1} a^{-3} w_{_2})^{2} w_{_1}^{-1},\quad w_{_4}=w_{_3}^{-2} w^{-4} w_{_3},\quad w_{_5}=c w_{_1} c,\quad
w_{_6}=w_{_5} w_{_3} w_{_5}^{-1},\]
 %w_{7} &=&  b^{4} a b^{-4} a^{-1} \,. \nonumber
 %w_{8} &=& w_{3} w_{4}^{-1} w_{2}^{-1} w_{1}^{-1}\,.\nonumber
 then the final three unipotent elements are
\[q_{_1}=w_{_4}^{-4} q_{_2},\quad q_{_2}=w_{_6} w_{_4}^{-32} w_{_3}^{-81},\quad
\quad
q_{_3}=w_{_3},\]which are respectively
{\tiny\[\begin{pmatrix}\begin {array}{rrrrr}1&0&0&16&0\\ \noalign{\medskip}0&1&0&0&
16\\ \noalign{\medskip}0&0&1&0&0\\ \noalign{\medskip}0&0&0&1&0
\\ \noalign{\medskip}0&0&0&0&1
\end {array}\end{pmatrix},\quad \begin{pmatrix}\begin {array}{rrrrr}1&0&-64&0&1024\\ \noalign{\medskip}0&1&0
&0&0\\ \noalign{\medskip}0&0&1&0&-32\\ \noalign{\medskip}0&0&0&1&0
\\ \noalign{\medskip}0&0&0&0&1
\end {array}\end{pmatrix},\quad \begin{pmatrix}\begin {array}{rrrrr}1&8&0&0&0\\ \noalign{\medskip}0&1&0&0&0
\\ \noalign{\medskip}0&0&1&0&0\\ \noalign{\medskip}0&0&0&1&8
\\ \noalign{\medskip}0&0&0&0&1
\end {array}\end{pmatrix}.\]}

\subsection{Arithmeticity of Case \ref{arithmetic-BS7}}\label{subsection-BS7}
In this case \[\alpha=\left(0, \frac{1}{3},\frac{1}{3}, \frac{2}{3}, \frac{2}{3}\right);\qquad \beta=\left(\frac{1}{2},\frac{1}{8},\frac{3}{8},\frac{5}{8},\frac{7}{8}\right)\]
\[f(x)=x^5+x^4+x^3-x^2-x-1;\quad g(x)=x^5+x^4+x+1.\]

The matrix (up to scalar multiplication) of the quadratic form $\mathrm{Q}$, 
with respect to the standard basis $\{e_1, e_2, e_3, e_4, e_5\}$, and  the change of basis matrix $X$, are
{\tiny\[\mathrm{Q}=\begin{pmatrix}\begin {array}{rrrrr}7&-5&-1&7&-9\\ \noalign{\medskip}-5&7&-5
&-1&7\\ \noalign{\medskip}-1&-5&7&-5&-1\\ \noalign{\medskip}7&-1&-5&7&
-5\\ \noalign{\medskip}-9&7&-1&-5&7
\end {array}
\end{pmatrix};\quad X=\begin{pmatrix}\begin {array}{rrrrr}1&0&0&0&1\\ \noalign{\medskip}2&0&4&-1&0
\\ \noalign{\medskip}1&1&4&0&0\\ \noalign{\medskip}0&2&4&0&-1
\\ \noalign{\medskip}0&1&0&1&0
\end {array}\end{pmatrix},\]}and the matrices of $\mathrm{Q}$, $A$ and $B$, with respect to the new basis, are $X^t\mathrm{Q}X$, $a=X^{-1}AX$ and $b=X^{-1}BX$; which are respectively
{\tiny\[\begin{pmatrix}\begin {array}{rrrrr}0&0&0&0&-4\\ \noalign{\medskip}0&0&0&-4&0
\\ \noalign{\medskip}0&0&-16&0&0\\ \noalign{\medskip}0&-4&0&0&0
\\ \noalign{\medskip}-4&0&0&0&0
\end {array}\end{pmatrix},\quad \begin{pmatrix}\begin {array}{rrrrr}-1&1&0&0&0\\ \noalign{\medskip}-1&0&0&0&0
\\ \noalign{\medskip}1&0&1&0&0\\ \noalign{\medskip}1&1&4&-1&-1
\\ \noalign{\medskip}1&0&0&1&0
\end {array}\end{pmatrix},\quad \begin{pmatrix}\begin {array}{rrrrr}-1&0&0&-1&0\\ \noalign{\medskip}-1&0&0&0
&0\\ \noalign{\medskip}1&0&1&0&0\\ \noalign{\medskip}1&1&4&-1&-1
\\ \noalign{\medskip}1&-1&0&0&0
\end{array}\end{pmatrix}.\]}

It is clear from the above computation that the $\Q$ - rank of the orthogonal group $\mathrm{SO}_\mathrm{Q}$ is two. If we denote by \[%c &=& a^{-1} b,\quad
d=(a b^{-4})^{2},\quad
%e & = & a^{2} b^{2} a^{-2} b^{-2},\quad
%u & = & a^{-1} b^{2} a b^{-2},\quad
%w & = & a^{-1} b^{4} a b^{-4},\quad
%r & = & (a b^{-4} )^{2},\quad
%t & = & w^{-1} u e^{-1},\quad
%z & = & b^4,\quad
w_{_1}=(a^{2} b^{4} a^{-1} d^{-1})^{-2},\quad
w_{_2}=(a^{-1} b a^{-4} b^{5})^{2},\quad
w_{_3}=[a^{-1}, b^{3}],\quad
w_{_4}=w_{_1} w_{_3}^{-1} w_{_1}^{-1}\]
\[w_{_5}=w_{_2} w_{_4}^{2},\quad
w_{_6}=w_{_3}^{4} w_{_5} w_{_3}^{4} w_{_5}^{-1},\quad
w_{_7}=[ w_{_3}^{4}, w_{_5} ],\]
 %w_{8} &=& w_{3} w_{4}^{-1} w_{2}^{-1} w_{1}^{-1}\,.\nonumber
then the final three unipotent elements are
\[q_{_1}=[ w_{_7}, w_{_5}^{-1}],\quad
q_{_2}=w_{_7}^{-2} q_{_1},\quad
q_{_3}=w_{_6}^{-2} q_{_2},\]which are respectively
{\tiny\[\begin{pmatrix}\begin {array}{rrrrr}1&0&0&-16&0\\ \noalign{\medskip}0&1&0&0&
16\\ \noalign{\medskip}0&0&1&0&0\\ \noalign{\medskip}0&0&0&1&0
\\ \noalign{\medskip}0&0&0&0&1
\end {array}\end{pmatrix},\quad \begin{pmatrix}\begin {array}{rrrrr}1&0&32&0&-128\\ \noalign{\medskip}0&1&0&0
&0\\ \noalign{\medskip}0&0&1&0&-8\\ \noalign{\medskip}0&0&0&1&0
\\ \noalign{\medskip}0&0&0&0&1
\end {array}\end{pmatrix},\quad \begin{pmatrix}\begin {array}{rrrrr}1&16&0&0&0\\ \noalign{\medskip}0&1&0&0&0
\\ \noalign{\medskip}0&0&1&0&0\\ \noalign{\medskip}0&0&0&1&-16
\\ \noalign{\medskip}0&0&0&0&1
\end {array}\end{pmatrix}.\]}

\subsection{Arithmeticity of Case \ref{arithmetic-BS8}}\label{subsection-BS8}
In this case \[\alpha=\left(0, \frac{1}{3},\frac{2}{3}, \frac{1}{4}, \frac{3}{4}\right);\qquad \beta=\left(\frac{1}{2},\frac{1}{8},\frac{3}{8},\frac{5}{8},\frac{7}{8}\right)\]
\[f(x)=x^5+x^3-x^2-1;\quad g(x)=x^5+x^4+x+1.\]

The matrix (up to scalar multiplication) of the quadratic form $\mathrm{Q}$, 
with respect to the standard basis $\{e_1, e_2, e_3, e_4, e_5\}$, and  the change of basis matrix $X$, are
{\tiny\[\mathrm{Q}=\begin{pmatrix}\begin {array}{rrrrr}3&-3&-1&5&-5\\ \noalign{\medskip}-3&3&-3
&-1&5\\ \noalign{\medskip}-1&-3&3&-3&-1\\ \noalign{\medskip}5&-1&-3&3&
-3\\ \noalign{\medskip}-5&5&-1&-3&3
\end {array}
\end{pmatrix};\quad X=\begin{pmatrix}\begin {array}{rrrrr}0&0&4/3&1&-2\\ \noalign{\medskip}0&6&4/3
&0&-1\\ \noalign{\medskip}2&6&0&-1&0\\ \noalign{\medskip}2&10&-8/9&-1/3&0\\ \noalign{\medskip}0&-6&0&1&-3
\end {array}\end{pmatrix},\]}and the matrices of $\mathrm{Q}$, $A$ and $B$, with respect to the new basis, are $X^t\mathrm{Q}X$, $a=X^{-1}AX$ and $b=X^{-1}BX$; which are respectively
{\tiny\[\begin{pmatrix}\begin {array}{rrrrr} 0&0&0&0&16\\ \noalign{\medskip}0&0&0&48&0
\\ \noalign{\medskip}0&0&-{{64}/{9}}&0&0\\ \noalign{\medskip}0&48
&0&0&0\\ \noalign{\medskip}16&0&0&0&0
\end {array}\end{pmatrix},\quad \begin{pmatrix}\begin {array}{rrrrr}-1/2&-7&-1/9&{{7}/{12}}&-11/4
\\ \noalign{\medskip}1/6&4/3&1/27&-{{7}/{36}}&{{19}/{36}}
\\ \noalign{\medskip}-3/2&-12&2/3&7/4&-{{17}/{4}}
\\ \noalign{\medskip}0&-6&-4/3&-1&5/3\\ \noalign{\medskip}-1&-8&-2/9&1
/6&-1/2
\end {array}\end{pmatrix},\quad \begin{pmatrix}\begin {array}{rrrrr}-1/2&-1&-1/9&-{{5}/{12}}&1/4
\\ \noalign{\medskip}1/6&1/3&1/27&-1/36&1/36\\ \noalign{\medskip}-3/2&
-3&2/3&1/4&1/4\\ \noalign{\medskip}0&-6&-4/3&-1&5/3
\\ \noalign{\medskip}-1&-8&-2/9&1/6&-1/2
\end{array}\end{pmatrix}.\]}

It is clear from the above computation that the $\Q$ - rank of the orthogonal group $\mathrm{SO}_\mathrm{Q}$ is two. If we denote by \[c=a^{-1} b,\quad
e=a b^{-1},\quad
w_{_1}=a^{2} b^{4} a^{-2},\quad
w_{_2}=w_{_1}^{-1} b^{4} w_{_1},\quad
w_{_3}=a^{-1} b^{4} a,\quad
w_{_4}=w_{_3} b^{6} w_{_3},\]
\[w_{_5}=[w_{_2}, w_{_4}^{2}],\quad w_{_6}=[c, w_{_4}^{-1}],\quad
w_{_7}=w_{_6}^{-16} w_{_5}^{3},\quad
w_{_8}=e w_{_6} e,\]
then the final three unipotent elements are
\[q_{_1}=w_{_5}^{-9} q_{_2},\quad 
q_{_2}=w_{_8}^{24} q_{_3},\quad q_{_3}=w_{_7} w_{_8}^{-8},\]which are respectively
{\tiny\[\begin{pmatrix}\begin {array}{rrrrr}1&0&0&144&0\\ \noalign{\medskip}0&1&0&0&
-48\\ \noalign{\medskip}0&0&1&0&0\\ \noalign{\medskip}0&0&0&1&0
\\ \noalign{\medskip}0&0&0&0&1
\end {array}\end{pmatrix},\quad \begin{pmatrix}\begin {array}{rrrrr}1&0&96&0&10368\\ \noalign{\medskip}0&1&0
&0&0\\ \noalign{\medskip}0&0&1&0&216\\ \noalign{\medskip}0&0&0&1&0
\\ \noalign{\medskip}0&0&0&0&1
\end {array}\end{pmatrix},\quad \begin{pmatrix}\begin {array}{rrrrr}1&-432&0&0&0\\ \noalign{\medskip}0&1&0&0
&0\\ \noalign{\medskip}0&0&1&0&0\\ \noalign{\medskip}0&0&0&1&144
\\ \noalign{\medskip}0&0&0&0&1
\end {array}\end{pmatrix}.\]}

\subsection{Arithmeticity of Case \ref{arithmetic-BS9} ($\Q$ - rank one case)}\label{subsection-BS9}
In this case \[\alpha=\left(0, \frac{1}{3},\frac{2}{3}, \frac{1}{4}, \frac{3}{4}\right);\qquad \beta=\left(\frac{1}{2},\frac{1}{12},\frac{5}{12},\frac{7}{12},\frac{11}{12}\right)\]
\[f(x)=x^5+x^3-x^2-1;\quad g(x)=x^5+x^4-x^3-x^2+x+1.\]

The matrix (up to scalar multiplication) of the quadratic form $\mathrm{Q}$, 
with respect to the standard basis $\{e_1, e_2, e_3, e_4, e_5\}$, and  the change of basis matrix $X$, are
{\tiny\[\mathrm{Q}=\begin{pmatrix}\begin {array}{rrrrr}1&-3&-1&3&-5\\ \noalign{\medskip}-3&1&-3
&-1&3\\ \noalign{\medskip}-1&-3&1&-3&-1\\ \noalign{\medskip}3&-1&-3&1&
-3\\ \noalign{\medskip}-5&3&-1&-3&1
\end {array}
\end{pmatrix};\quad X=\begin{pmatrix}\begin {array}{rrrrr}-1&0&0&0&-1\\ \noalign{\medskip}0&-1&-2&
1&-1\\ \noalign{\medskip}-1&-2&-1&-1&1\\ \noalign{\medskip}0&-4&1&-1&2
\\ \noalign{\medskip}0&1&2&-1&-1
\end {array}\end{pmatrix},\]}and the matrices of $\mathrm{Q}$, $A$ and $B$, with respect to the new basis, are $X^t\mathrm{Q}X$, $a=X^{-1}AX$ and $b=X^{-1}BX$; which are respectively
{\tiny\[\begin{pmatrix}\begin {array}{rrrrr} 0&0&0&0&-12\\ \noalign{\medskip}0&-24&-
12&-12&0\\ \noalign{\medskip}0&-12&-24&12&0\\ \noalign{\medskip}0&-12&
12&-8&0\\ \noalign{\medskip}-12&0&0&0&0
\end {array}\end{pmatrix},\quad \begin{pmatrix}\begin {array}{rrrrr}-1/2&-3&-3/2&1/2&3/2
\\ \noalign{\medskip}1/2&1&1/2&1/6&-1/2\\ \noalign{\medskip}0&0&0&-1/3
&1\\ \noalign{\medskip}0&3&0&0&0\\ \noalign{\medskip}1/2&2&-1/2&1/2&-1
/2
\end {array}\end{pmatrix},\quad \begin{pmatrix}\begin {array}{rrrrr}-1/2&-2&1/2&-1/2&1/2
\\ \noalign{\medskip}1/2&1&1/2&1/6&-1/2\\ \noalign{\medskip}0&0&0&-1/3
&1\\ \noalign{\medskip}0&3&0&0&0\\ \noalign{\medskip}1/2&3&3/2&-1/2&-3
/2
\end{array}\end{pmatrix}.\]}

A computation shows that the $\Q$ - rank of the orthogonal group $\mathrm{SO}_\mathrm{Q}$ is one. If we denote by \[c=a^{-1} b,\quad
d=aba,\quad
e=[a^{-3}, b^{-1}],\quad
w_{_1}=c d c a^{-3},\quad
w_{_2}=a^{3} e^{-1} a^{-3} d c a^{-3},\quad
w_{_3}=w_{_1}^{-2} w_{_2}^{3},\]
\[w_{_4}=w_{_3} w_{_2}^{-1},\quad w_{_5}=w_{_3}^{-1} w_{_4}^{3},\quad
w_{_6}=[e, w_{_1}^{2}],\quad
w_{_7}=w_{_5}^{3} w_{_6},\]
%w_{8} &=&  e w_{6} e\,.\nonumber
then the final three unipotent elements are
\[q_{_1}=w_{_1}^{8} q_{_2},\quad
q_{_2}=w_{_1}^{4} w_{_7},\quad
q_{_3}=w_{_3}^{4} q_{_1}^{-1},\]which are respectively
{\tiny\[\begin{pmatrix}\begin {array}{rrrrr}1&0&0&24&216\\ \noalign{\medskip}0&1&0&0
&-18\\ \noalign{\medskip}0&0&1&0&18\\ \noalign{\medskip}0&0&0&1&18
\\ \noalign{\medskip}0&0&0&0&1
\end {array}\end{pmatrix},\quad \begin{pmatrix}\begin {array}{rrrrr}1&0&24&0&24\\ \noalign{\medskip}0&1&0&0&
-10\\ \noalign{\medskip}0&0&1&0&2\\ \noalign{\medskip}0&0&0&1&18
\\ \noalign{\medskip}0&0&0&0&1
\end {array}\end{pmatrix},\quad \begin{pmatrix}\begin {array}{rrrrr}1&24&0&0&24\\ \noalign{\medskip}0&1&0&0&
2\\ \noalign{\medskip}0&0&1&0&-10\\ \noalign{\medskip}0&0&0&1&-18
\\ \noalign{\medskip}0&0&0&0&1
\end {array}\end{pmatrix}.\]}

\subsection{Arithmeticity of Case \ref{arithmetic-BS10}}\label{subsection-BS10}
In this case \[\alpha=\left(0, \frac{1}{3},\frac{2}{3}, \frac{1}{6}, \frac{5}{6}\right);\qquad \beta=\left(\frac{1}{2},\frac{1}{4},\frac{1}{4},\frac{3}{4},\frac{3}{4}\right)\]
\[f(x)=x^5-x^4+x^3-x^2+x-1;\quad g(x)=x^5+x^4+2x^3+2x^2+x+1.\]

The matrix (up to scalar multiplication) of the quadratic form $\mathrm{Q}$, 
with respect to the standard basis $\{e_1, e_2, e_3, e_4, e_5\}$, and  the change of basis matrix $X$, are
{\tiny\[\mathrm{Q}=\begin{pmatrix}\begin {array}{rrrrr}1&-7&1&17&1\\ \noalign{\medskip}-7&1&-7&
1&17\\ \noalign{\medskip}1&-7&1&-7&1\\ \noalign{\medskip}17&1&-7&1&-7
\\ \noalign{\medskip}1&17&1&-7&1
\end {array}
\end{pmatrix};\quad X=\begin{pmatrix}\begin {array}{rrrrr}0&1&0&-1&2\\ \noalign{\medskip}1&1&0&-1&
1\\ \noalign{\medskip}0&2&0&-1&4\\ \noalign{\medskip}-1&3&-24&4&0
\\ \noalign{\medskip}0&1&0&-1&1
\end {array}\end{pmatrix},\]}and the matrices of $\mathrm{Q}$, $A$ and $B$, with respect to the new basis, are $X^t\mathrm{Q}X$, $a=X^{-1}AX$ and $b=X^{-1}BX$; which are respectively
{\tiny\[\begin{pmatrix}\begin {array}{rrrrr} 0&0&0&0&-24\\ \noalign{\medskip}0&0&0&-
24&0\\ \noalign{\medskip}0&0&576&0&0\\ \noalign{\medskip}0&-24&0&0&0
\\ \noalign{\medskip}-24&0&0&0&0
\end {array}\end{pmatrix},\quad \begin{pmatrix}\begin {array}{rrrrr}1&-4&24&-3&0\\ \noalign{\medskip}-1&7&-
48&7&1\\ \noalign{\medskip}0&1&-7&1&0\\ \noalign{\medskip}1&0&0&0&0
\\ \noalign{\medskip}1&-3&24&-4&0
\end {array}\end{pmatrix},\quad \begin{pmatrix}\begin {array}{rrrrr}1&-2&24&-5&2\\ \noalign{\medskip}-1&6&-
48&8&0\\ \noalign{\medskip}0&1&-7&1&0\\ \noalign{\medskip}1&1&0&-1&1
\\ \noalign{\medskip}1&-3&24&-4&0
\end{array}\end{pmatrix}.\]}

It is clear from the above computation that the $\Q$ - rank of the orthogonal group $\mathrm{SO}_\mathrm{Q}$ is two. If we denote by \[d=[ a , b^{-1} ],\quad
e=[ b^{-1} , a^{2} ],\quad
r=a b a^{-1},\quad
s=b a^{2} b^{-1},\quad
t=[ s,  d^{-1} ],\quad
w_{_1}=b^{-1} a^{-1} b,\]
\[w_{_2}={[e^{-1} , t^{-1}]}^{2},\quad
w_{_3}=d^{-12} w_{_2},\quad
w_{_4}=s w_{_2} s^{-1},\quad
w_{_5}=d^{2} w_{_1}^{2} d^{-1},\quad
w_{_6}=w_{_4} w_{_5}^{2} w_{_4}^{-1} w_{_5}^{2},\]
%w_{7} & = & w_{5}^{3} w_{6}\,. \nonumber
%w_{8} &=&  e w_{6} e\,.\nonumber
then the final three unipotent elements are
\[q_{_1}=w_{_3},\quad q_{_2}=[ w_{_6}, w_{_2} ] w_{_3}^{-1728},\quad
q_{_3}=w_{_2},\]which are respectively
{\tiny\[\begin{pmatrix}\begin {array}{rrrrr}1&0&0&-12&0\\ \noalign{\medskip}0&1&0&0&
12\\ \noalign{\medskip}0&0&1&0&0\\ \noalign{\medskip}0&0&0&1&0
\\ \noalign{\medskip}0&0&0&0&1
\end {array}\end{pmatrix},\quad \begin{pmatrix}\begin {array}{rrrrr}1&0&-3456&0&248832\\ \noalign{\medskip}0
&1&0&0&0\\ \noalign{\medskip}0&0&1&0&-144\\ \noalign{\medskip}0&0&0&1&0
\\ \noalign{\medskip}0&0&0&0&1
\end {array}\end{pmatrix},\quad \begin{pmatrix}\begin {array}{rrrrr}1&-12&0&0&0\\ \noalign{\medskip}0&1&0&0&0
\\ \noalign{\medskip}0&0&1&0&0\\ \noalign{\medskip}0&0&0&1&12
\\ \noalign{\medskip}0&0&0&0&1
\end {array}\end{pmatrix}.\]}

\subsection{Arithmeticity of Case \ref{arithmetic-BS11}}\label{subsection-BS11}
In this case \[\alpha=\left(0, \frac{1}{3},\frac{2}{3}, \frac{1}{6}, \frac{5}{6}\right);\qquad \beta=\left(\frac{1}{2},\frac{1}{8},\frac{3}{8},\frac{5}{8},\frac{7}{8}\right)\]
\[f(x)=x^5-x^4+x^3-x^2+x-1;\quad g(x)=x^5+x^4+x+1.\]

The matrix (up to scalar multiplication) of the quadratic form $\mathrm{Q}$, 
with respect to the standard basis $\{e_1, e_2, e_3, e_4, e_5\}$, and  the change of basis matrix $X$, are
{\tiny\[\mathrm{Q}=\begin{pmatrix}\begin {array}{rrrrr}1&5&1&-7&1\\ \noalign{\medskip}5&1&5&1&-
7\\ \noalign{\medskip}1&5&1&5&1\\ \noalign{\medskip}-7&1&5&1&5
\\ \noalign{\medskip}1&-7&1&5&1
\end {array}
\end{pmatrix};\quad X=\begin{pmatrix}\begin {array}{rrrrr}0&2&0&-1&2\\ \noalign{\medskip}-1&2&0&-1
&1\\ \noalign{\medskip}0&2&0&0&0\\ \noalign{\medskip}1&-8&-12&-1&-2
\\ \noalign{\medskip}0&2&0&-1&3
\end {array}\end{pmatrix},\]}and the matrices of $\mathrm{Q}$, $A$ and $B$, with respect to the new basis, are $X^t\mathrm{Q}X$, $a=X^{-1}AX$ and $b=X^{-1}BX$; which are respectively
{\tiny\[\begin{pmatrix}\begin {array}{rrrrr} 0&0&0&0&12\\ \noalign{\medskip}0&0&0&-24
&0\\ \noalign{\medskip}0&0&144&0&0\\ \noalign{\medskip}0&-24&0&0&0
\\ \noalign{\medskip}12&0&0&0&0
\end {array}\end{pmatrix},\quad \begin{pmatrix}\begin {array}{rrrrr}-1&10&12&0&6\\ \noalign{\medskip}-1/2&2&0
&-1&2\\ \noalign{\medskip}0&2&5&1&0\\ \noalign{\medskip}1&-14&-24&-3&-
3\\ \noalign{\medskip}1&-8&-12&-1&-2
\end {array}\end{pmatrix},\quad \begin{pmatrix}\begin {array}{rrrrr}-1&6&12&2&0\\ \noalign{\medskip}-1/2&1&0
&-1/2&1/2\\ \noalign{\medskip}0&2&5&1&0\\ \noalign{\medskip}1&-12&-24&
-4&0\\ \noalign{\medskip}1&-8&-12&-1&-2
\end{array}\end{pmatrix}.\]}

It is clear from the above computation that the $\Q$ - rank of the orthogonal group $\mathrm{SO}_\mathrm{Q}$ is two. If we denote by \[d={[ b^{-1}, a ]}^{2},\quad
w_{_1}=a b^{4} a^{-1},\quad
w_{_2}=[ a^{2},  w_{_1}^{-1} ],\quad
w_{_3}=w_{_2} d w_{_2}^{-1} d,\quad
w_{_4}=w_{_3} d^{-2},\]
then the final three unipotent elements are
\[q_{_1}=w_{_4}^{2} q_{_2}^{-1},\quad q_{_2}=w_{_2}^{-1} w_{_4} w_{_2} w_{_4},\quad
q_{_3}=(q_{_1} d^{6})^{2} q_{_1}^{-1},\]which are respectively
{\tiny\[\begin{pmatrix}\begin {array}{rrrrr}1&0&0&24&0\\ \noalign{\medskip}0&1&0&0&
12\\ \noalign{\medskip}0&0&1&0&0\\ \noalign{\medskip}0&0&0&1&0
\\ \noalign{\medskip}0&0&0&0&1
\end {array}\end{pmatrix},\quad \begin{pmatrix}\begin {array}{rrrrr}1&0&-48&0&-96\\ \noalign{\medskip}0&1&0&0
&0\\ \noalign{\medskip}0&0&1&0&4\\ \noalign{\medskip}0&0&0&1&0
\\ \noalign{\medskip}0&0&0&0&1
\end {array}\end{pmatrix},\quad \begin{pmatrix}\begin {array}{rrrrr}1&48&0&0&0\\ \noalign{\medskip}0&1&0&0&0
\\ \noalign{\medskip}0&0&1&0&0\\ \noalign{\medskip}0&0&0&1&24
\\ \noalign{\medskip}0&0&0&0&1
\end {array}\end{pmatrix}.\]}

\subsection{Arithmeticity of Case \ref{arithmetic-BS12}}\label{subsection-BS12}
In this case \[\alpha=\left(0, \frac{1}{3},\frac{2}{3}, \frac{1}{6}, \frac{5}{6}\right);\qquad \beta=\left(\frac{1}{2},\frac{1}{12},\frac{5}{12},\frac{7}{12},\frac{11}{12}\right)\]
\[f(x)=x^5-x^4+x^3-x^2+x-1;\quad g(x)=x^5 +x^4-x^3-x^2+x+1.\]

The matrix (up to scalar multiplication) of the quadratic form $\mathrm{Q}$, 
with respect to the standard basis $\{e_1, e_2, e_3, e_4, e_5\}$, and  the change of basis matrix $X$, are
{\tiny\[\mathrm{Q}=\begin{pmatrix}\begin {array}{rrrrr}1&2&1&-1&1\\ \noalign{\medskip}2&1&2&1&-
1\\ \noalign{\medskip}1&2&1&2&1\\ \noalign{\medskip}-1&1&2&1&2
\\ \noalign{\medskip}1&-1&1&2&1
\end {array}
\end{pmatrix};\quad X=\begin{pmatrix}\begin {array}{rrrrr}0&1&0&-1/2&1\\ \noalign{\medskip}-1&1&0&
-1/2&1/2\\ \noalign{\medskip}0&1&0&1/2&-1\\ \noalign{\medskip}1&-4&-3&0
&-3/2\\ \noalign{\medskip}0&1&0&-1/2&2
\end {array}\end{pmatrix},\]}and the matrices of $\mathrm{Q}$, $A$ and $B$, with respect to the new basis, are $X^t\mathrm{Q}X$, $a=X^{-1}AX$ and $b=X^{-1}BX$; which are respectively
{\tiny\[\begin{pmatrix}\begin {array}{rrrrr} 0&0&0&0&3\\ \noalign{\medskip}0&0&0&-3&0
\\ \noalign{\medskip}0&0&9&0&0\\ \noalign{\medskip}0&-3&0&0&0
\\ \noalign{\medskip}3&0&0&0&0
\end {array}\end{pmatrix},\quad \begin{pmatrix}\begin {array}{rrrrr}-1/2&3&3/2&-1/2&{{15}/{4}}
\\ \noalign{\medskip}-1/2&3/2&0&-3/4&9/4\\ \noalign{\medskip}0&1&2&1/2
&0\\ \noalign{\medskip}1&-7&-6&-1/2&-5/2\\ \noalign{\medskip}1&-4&-3&0
&-3/2
\end {array}\end{pmatrix},\quad \begin{pmatrix}\begin {array}{rrrrr}-1/2&1&3/2&1/2&-1/4\\ \noalign{\medskip}
-1/2&1/2&0&-1/4&1/4\\ \noalign{\medskip}0&1&2&1/2&0
\\ \noalign{\medskip}1&-5&-6&-3/2&3/2\\ \noalign{\medskip}1&-4&-3&0&-3
/2
\end{array}\end{pmatrix}.\]}

It is clear from the above computation that the $\Q$ - rank of the orthogonal group $\mathrm{SO}_\mathrm{Q}$ is two. If we denote by \[d=[ b^{-1}, a ],\quad
e=a b^{6} a^{-1},\quad
w_{_1}=[ a^{2},  e^{-1} ],\quad
w_{_2}=w_{_1} d w_{_1}^{-1} d,\]
then the final three unipotent elements are
\[q_{_1} =  \bigg(\big( q_{_2}^{-1} d q_{_2}\big)^{2} w_{_2}^{-1}\bigg)^{2} q_{_2},\quad q_{_2}=w_{_1}^{-1} w_{_2} d^{-2} w_{_1} w_{_2} d^{-2},\quad
q_{_3}  = (q_{_2}^{-1} d q_{_2})^{-96} q_{_1},\]
which are respectively
{\tiny\[\begin{pmatrix}\begin {array}{rrrrr}1&0&0&-96&0\\ \noalign{\medskip}0&1&0&0&
-96\\ \noalign{\medskip}0&0&1&0&0\\ \noalign{\medskip}0&0&0&1&0
\\ \noalign{\medskip}0&0&0&0&1
\end {array}\end{pmatrix},\quad \begin{pmatrix}\begin {array}{rrrrr}1&0&-48&0&-384\\ \noalign{\medskip}0&1&0
&0&0\\ \noalign{\medskip}0&0&1&0&16\\ \noalign{\medskip}0&0&0&1&0
\\ \noalign{\medskip}0&0&0&0&1
\end {array}\end{pmatrix},\quad \begin{pmatrix}\begin {array}{rrrrr}1&-192&0&0&0\\ \noalign{\medskip}0&1&0&0
&0\\ \noalign{\medskip}0&0&1&0&0\\ \noalign{\medskip}0&0&0&1&-192
\\ \noalign{\medskip}0&0&0&0&1
\end {array}\end{pmatrix}.\]}

\subsection{Arithmeticity of Case \ref{arithmetic-BS13}}\label{subsection-BS13}
In this case \[\alpha=\left(0, \frac{1}{5},\frac{2}{5}, \frac{3}{5}, \frac{4}{5}\right);\qquad \beta=\left(\frac{1}{2},\frac{1}{3},\frac{1}{3},\frac{2}{3},\frac{2}{3}\right)\]
\[f(x)= x^5-1 ;\quad g(x)=x^5+3x^4+5x^3+5x^2+3x+1.\]

The matrix (up to scalar multiplication) of the quadratic form $\mathrm{Q}$, 
with respect to the standard basis $\{e_1, e_2, e_3, e_4, e_5\}$, and  the change of basis matrix $X$, are
{\tiny\[\mathrm{Q}=\begin{pmatrix}\begin {array}{rrrrr}29&-25&11&11&-25\\ \noalign{\medskip}-25
&29&-25&11&11\\ \noalign{\medskip}11&-25&29&-25&11
\\ \noalign{\medskip}11&11&-25&29&-25\\ \noalign{\medskip}-25&11&11&-
25&29
\end {array}
\end{pmatrix};\quad X=\begin{pmatrix}\begin {array}{rrrrr}-1&-2&2&3&6\\ \noalign{\medskip}-1&-1&2&
4&4\\ \noalign{\medskip}1&0&0&0&1\\ \noalign{\medskip}1&-1&0&-1&2
\\ \noalign{\medskip}0&-2&2&0&5
\end {array}\end{pmatrix},\]}and the matrices of $\mathrm{Q}$, $A$ and $B$, with respect to the new basis, are $X^t\mathrm{Q}X$, $a=X^{-1}AX$ and $b=X^{-1}BX$; which are respectively
{\tiny\[\begin{pmatrix}\begin {array}{rrrrr} 0&0&0&0&18\\ \noalign{\medskip}0&0&0&-18
&0\\ \noalign{\medskip}0&0&36&0&0\\ \noalign{\medskip}0&-18&0&0&0
\\ \noalign{\medskip}18&0&0&0&0
\end {array}\end{pmatrix},\quad \begin{pmatrix}\begin {array}{rrrrr}0&0&0&3&2\\ \noalign{\medskip}-3&-2&4&4&
4\\ \noalign{\medskip}0&0&-1&1&0\\ \noalign{\medskip}0&0&0&1&1
\\ \noalign{\medskip}-1&-1&2&1&2
\end {array}\end{pmatrix},\quad \begin{pmatrix}\begin {array}{rrrrr}0&12&-12&3&-28\\ \noalign{\medskip}-3&-8
&10&4&19\\ \noalign{\medskip}0&2&-3&1&-5\\ \noalign{\medskip}0&4&-4&1&
-9\\ \noalign{\medskip}-1&-3&4&1&7
\end{array}\end{pmatrix}.\]}

It is clear from the above computation that the $\Q$ - rank of the orthogonal group $\mathrm{SO}_\mathrm{Q}$ is two. If we denote by \[%c = a^{-1} b,\quad
w_{_1}=  [ a^{-1}, b^{-6} ],\quad
w_{_2} = [ a^{-1}, b^{3}],\quad
w_{_3} = a^{-3} b^{-6} a^{3},\quad
w_{_4} = w_{_2}^{2} w_{1}^{-1} w_{_3} w_{_2}^{-2} w_{_1} w_{_3}^{-1},\]
\[w_{_5} = w_{_4}^{-2}  w_{_1} w_{_4}^{-1} w_{_2}^{-2},\]
then the final three unipotent elements are
\[q_{_1}= w_{_3},\quad
q_{_2}  = w_{_5}^{3} w_{_4}^{10},\quad
q_{_3} = w_{_4},\]
which are respectively
{\tiny\[\begin{pmatrix}\begin {array}{rrrrr}1&0&0&6&0\\ \noalign{\medskip}0&1&0&0&6
\\ \noalign{\medskip}0&0&1&0&0\\ \noalign{\medskip}0&0&0&1&0
\\ \noalign{\medskip}0&0&0&0&1
\end {array}\end{pmatrix},\quad \begin{pmatrix}\begin {array}{rrrrr}1&0&24&0&-144\\ \noalign{\medskip}0&1&0&0
&0\\ \noalign{\medskip}0&0&1&0&-12\\ \noalign{\medskip}0&0&0&1&0
\\ \noalign{\medskip}0&0&0&0&1
\end {array}\end{pmatrix},\quad \begin{pmatrix}\begin {array}{rrrrr}1&-48&0&0&0\\ \noalign{\medskip}0&1&0&0&0
\\ \noalign{\medskip}0&0&1&0&0\\ \noalign{\medskip}0&0&0&1&-48
\\ \noalign{\medskip}0&0&0&0&1
\end {array}\end{pmatrix}.\]}

\subsection{Arithmeticity of Case \ref{arithmetic-BS14} ($\Q$ - rank one case)}\label{subsection-BS14}
In this case \[\alpha=\left(0, \frac{1}{5},\frac{2}{5}, \frac{3}{5}, \frac{4}{5}\right);\qquad \beta=\left(\frac{1}{2},\frac{1}{3},\frac{2}{3},\frac{1}{4},\frac{3}{4}\right)\]
\[f(x)= x^5-1 ;\quad g(x)=x^5 +2x^4+3x^3+3x^2+2x+1.\]

The matrix (up to scalar multiplication) of the quadratic form $\mathrm{Q}$, 
with respect to the standard basis $\{e_1, e_2, e_3, e_4, e_5\}$, and  the change of basis matrix $X$, are
{\tiny\[\mathrm{Q}=\begin{pmatrix}\begin {array}{rrrrr}5&-7&5&5&-7\\ \noalign{\medskip}-7&5&-7&
5&5\\ \noalign{\medskip}5&-7&5&-7&5\\ \noalign{\medskip}5&5&-7&5&-7
\\ \noalign{\medskip}-7&5&5&-7&5
\end {array}
\end{pmatrix};\quad X=\begin{pmatrix}\begin {array}{rrrrr}0&-2&3&-1&2\\ \noalign{\medskip}1&-4&6&-
2&3\\ \noalign{\medskip}0&-3&7&-3&4\\ \noalign{\medskip}-1&-1&5&-5&2
\\ \noalign{\medskip}0&-2&3&-1&1
\end {array}\end{pmatrix},\]}and the matrices of $\mathrm{Q}$, $A$ and $B$, with respect to the new basis, are $X^t\mathrm{Q}X$, $a=X^{-1}AX$ and $b=X^{-1}BX$; which are respectively
{\tiny\[\begin{pmatrix}\begin {array}{rrrrr} 0&0&0&0&-12\\ \noalign{\medskip}0&24&-36
&12&0\\ \noalign{\medskip}0&-36&24&0&0\\ \noalign{\medskip}0&12&0&24&0
\\ \noalign{\medskip}-12&0&0&0&0
\end {array}\end{pmatrix},\quad \begin{pmatrix}\begin {array}{rrrrr}1&1&-5&5&-1\\ \noalign{\medskip}1&0&-1&1
&0\\ \noalign{\medskip}0&0&2&-3&1\\ \noalign{\medskip}0&0&1&-2&0
\\ \noalign{\medskip}1&-1&-2&4&-1
\end {array}\end{pmatrix},\quad \begin{pmatrix}\begin {array}{rrrrr}1&-3&1&3&1\\ \noalign{\medskip}1&-2&2&0&
1\\ \noalign{\medskip}0&0&2&-3&1\\ \noalign{\medskip}0&0&1&-2&0
\\ \noalign{\medskip}1&-1&-2&4&-1
\end{array}\end{pmatrix}.\]}

A computation shows that the $\Q$ - rank of the orthogonal group $\mathrm{SO}_\mathrm{Q}$ is one. If we denote by \[c = a^{-1} b,\quad
d = [a, b^{-1}],\quad
e = [a^{-4}, b^{-4}],\quad
w_{_1} = d (e d e^{-1})^{2},\quad
w_{_2} = c [a^{2}, b],\]
\[w_{_3} = w_{_2}^{2} (d^{2} w_{_2}^2)^{-1} w_{_2}^{2},\quad
w_{_4} = w_{_1}^{-4} w_{_3}^{3},\quad
w_{_5} = \big(w_{_2}  w_{_3}^{-5} w_{_2}^{-1}\big)^{2} w_{4}^{-2},\quad
w_{_6} = d^{98} w_{_5},\]
then the final three unipotent elements are
\[q_{_1}= w_{_6}^{-432} q_{_2}^{294},\quad
q_{_2} = w_{_1}^{46} w_{6},\quad
q_{_3}  = w_{_4}^{432} q_{_2}^{30},\]which are respectively
{\tiny\[\begin{pmatrix}\begin {array}{rrrrr}1&0&0&-59616&740430720
\\ \noalign{\medskip}0&1&0&0&-9936\\ \noalign{\medskip}0&0&1&0&-14904
\\ \noalign{\medskip}0&0&0&1&-24840\\ \noalign{\medskip}0&0&0&0&1
\end {array}\end{pmatrix},\quad \begin{pmatrix}\begin {array}{rrrrr}1&0&-432&0&-23328\\ \noalign{\medskip}0&
1&0&0&216\\ \noalign{\medskip}0&0&1&0&108\\ \noalign{\medskip}0&0&0&1&
-108\\ \noalign{\medskip}0&0&0&0&1
\end {array}\end{pmatrix},\quad \begin{pmatrix}\begin {array}{rrrrr}1&1296&0&0&-279936\\ \noalign{\medskip}0
&1&0&0&-432\\ \noalign{\medskip}0&0&1&0&-648\\ \noalign{\medskip}0&0&0
&1&216\\ \noalign{\medskip}0&0&0&0&1
\end {array}\end{pmatrix}.\]}

\subsection{Arithmeticity of Case \ref{arithmetic-BS15} [(f(x), g(x)), (-f(-x), -g(-x)) case]}\label{subsection-BS15}
In this case \[\alpha=\left(0, \frac{1}{6},\frac{1}{6}, \frac{5}{6}, \frac{5}{6}\right);\qquad \beta=\left(\frac{1}{2},\frac{1}{3},\frac{1}{3},\frac{2}{3},\frac{2}{3}\right)\]
\[f(x)= x^5-3x^4+5x^3-5x^2+3x-1;\quad g(x)=x^5 +3x^4+5x^3+5x^2+3x+1.\]

The matrix (up to scalar multiplication) of the quadratic form $\mathrm{Q}$, 
with respect to the standard basis $\{e_1, e_2, e_3, e_4, e_5\}$, and  the change of basis matrix $X$, are
{\tiny\[\mathrm{Q}=\begin{pmatrix}\begin {array}{rrrrr}4&0&-5&0&13\\ \noalign{\medskip}0&4&0&-5
&0\\ \noalign{\medskip}-5&0&4&0&-5\\ \noalign{\medskip}0&-5&0&4&0
\\ \noalign{\medskip}13&0&-5&0&4
\end {array}
\end{pmatrix};\quad X=\begin{pmatrix}\begin {array}{rrrrr}-1&-5&1/4&7/4&-5\\ \noalign{\medskip}3&
13&-3/4&-19/4&15\\ \noalign{\medskip}-4&-18&1&21/2&-26
\\ \noalign{\medskip}3&14&-1/2&-13/2&18\\ \noalign{\medskip}-1&-4&1/2&
8&-20
\end {array}\end{pmatrix},\]}and the matrices of $\mathrm{Q}$, $A$ and $B$, with respect to the new basis, are $X^t\mathrm{Q}X$, $a=X^{-1}AX$ and $b=X^{-1}BX$; which are respectively
{\tiny\[\begin{pmatrix}\begin {array}{rrrrr} 0&0&0&0&-18\\ \noalign{\medskip}0&0&0&-
36&0\\ \noalign{\medskip}0&0&1/2&0&0\\ \noalign{\medskip}0&-36&0&0&0
\\ \noalign{\medskip}-18&0&0&0&0
\end {array}\end{pmatrix},\qquad \begin{pmatrix}\begin {array}{rrrrr}-5/2&-23/2&5/8&-{{15}/{8}}&3/2
\\ \noalign{\medskip}{{9}/{16}}&{{41}/{16}}&-{{11}/{64}
}&-{{63}/{64}}&{{51}/{16}}\\ \noalign{\medskip}-9/2&-{{45}/{2}}&-1/8&-{{117}/{8}}&{{81}/{2}}\\ \noalign{\medskip}1/
4&1/4&-3/16&-{{7}/{16}}&11/4\\ \noalign{\medskip}0&-1/2&-1/8&-{{9}/{8}}&7/2
\end {array}\end{pmatrix},\]\[\begin{pmatrix}\begin {array}{rrrrr}15/2&{{57}/{2}}&-{{35}/{8}}&-{{655}/{8}}&{{403}/{2}}\\ \noalign{\medskip}-{{19}/{16}}
&-{{71}/{16}}&{{45}/{64}}&{{833}/{64}}&-{{509}/{16
}}\\ \noalign{\medskip}{{27}/{2}}&{{99}/{2}}&-{{73}/{8}
}&-{{1269}/{8}}&{{801}/{2}}\\ \noalign{\medskip}{{21}/{
4}}&{{81}/{4}}&-{{43}/{16}}&-{{647}/{16}}&{{411}/{
4}}\\ \noalign{\medskip}2&15/2&-{{9}/{8}}&-{{137}/{8}}&{{87}/{2}}
\end{array}\end{pmatrix}.\]}

It is clear from the above computation that the $\Q$ - rank of the orthogonal group $\mathrm{SO}_\mathrm{Q}$ is two. If we denote by \[%c = a^{-1} b,\quad
d = a b a^{-2},\quad
e = a b a,\quad
t = b a^{-1} b,\quad
w_{_1} = [ a^{-6}, e],\quad
w_{_2} = a^{-48} w_{_1},\quad
w_{_3} = t b^{-6} t^{-1},\]
\[w_{_4} = d w_{_3} d w_{_3}^{-1},\quad
w_{_5} = w_{_3} d w_{_3}^{-1} d,\quad
w_{_6} = [ w_{_5}, a^{6} ],\quad
w_{_7} = [ w_{_4}, a^{6} ],\]
then the final three unipotent elements are
\[q_{_1}= w_{_6} w_{_7}^{-1} w_{_1}^{-10} w_{_2} w_{_1}^{-2} w_{_2}^{3},\quad
q_{_2} = d q_{_1} d,\quad
q_{_3}  = w_{_2}^{384} q_{_1}^{-60},\]
which are respectively
{\tiny\[\begin{pmatrix}\begin {array}{rrrrr}1&0&0&-768&0\\ \noalign{\medskip}0&1&0&0
&384\\ \noalign{\medskip}0&0&1&0&0\\ \noalign{\medskip}0&0&0&1&0
\\ \noalign{\medskip}0&0&0&0&1
\end {array}\end{pmatrix},\quad \begin{pmatrix}\begin {array}{rrrrr}1&0&-4608&0&382205952
\\ \noalign{\medskip}0&1&0&0&0\\ \noalign{\medskip}0&0&1&0&-165888
\\ \noalign{\medskip}0&0&0&1&0\\ \noalign{\medskip}0&0&0&0&1
\end {array}\end{pmatrix},\quad \begin{pmatrix}\begin {array}{rrrrr}1&-3072&0&0&0\\ \noalign{\medskip}0&1&0&0
&0\\ \noalign{\medskip}0&0&1&0&0\\ \noalign{\medskip}0&0&0&1&1536
\\ \noalign{\medskip}0&0&0&0&1
\end {array}\end{pmatrix}.\]}

\subsection{Arithmeticity of Case \ref{arithmetic-BS16}}\label{subsection-BS16}
In this case \[\alpha=\left(0, \frac{1}{6},\frac{1}{6}, \frac{5}{6}, \frac{5}{6}\right);\qquad \beta=\left(\frac{1}{2},\frac{1}{3},\frac{2}{3},\frac{1}{4},\frac{3}{4}\right)\]
\[f(x)= x^5-3x^4+5x^3-5x^2+3x-1;\quad g(x)=x^5 +2x^4+3x^3+3x^2+2x+1.\]

The matrix (up to scalar multiplication) of the quadratic form $\mathrm{Q}$, 
with respect to the standard basis $\{e_1, e_2, e_3, e_4, e_5\}$, and  the change of basis matrix $X$, are
{\tiny\[\mathrm{Q}=\begin{pmatrix}\begin {array}{rrrrr}13&2&-14&-7&22\\ \noalign{\medskip}2&13&
2&-14&-7\\ \noalign{\medskip}-14&2&13&2&-14\\ \noalign{\medskip}-7&-14
&2&13&2\\ \noalign{\medskip}22&-7&-14&2&13
\end {array}
\end{pmatrix};\quad X=\begin{pmatrix}\begin {array}{rrrrr}1&7&0&19&36\\ \noalign{\medskip}-1&-8&3/49&{{166}/{7}}&42\\ \noalign{\medskip}0&2&1/49&{
{242}/{7}}&64\\ \noalign{\medskip}1&5&-1/49&{{185}/{7}}&48
\\ \noalign{\medskip}-1&-6&{{4}/{49}}&{{114}/{7}}&26
\end {array}\end{pmatrix},\]}and the matrices of $\mathrm{Q}$, $A$ and $B$, with respect to the new basis, are $X^t\mathrm{Q}X$, $a=X^{-1}AX$ and $b=X^{-1}BX$; which are respectively
{\tiny\[\begin{pmatrix}\begin {array}{rrrrr} 0&0&0&0&-72\\ \noalign{\medskip}0&0&0&
252&0\\ \noalign{\medskip}0&0&{{3}/{49}}&0&0\\ \noalign{\medskip}0
&252&0&0&0\\ \noalign{\medskip}-72&0&0&0&0
\end {array}\end{pmatrix},\qquad \begin{pmatrix}\begin {array}{rrrrr}20&127&-{{73}/{49}}&-{{2042}/{
7}}&-464\\ \noalign{\medskip}-3&-19&{{11}/{49}}&{{306}/{7}}&
{{3392}/{49}}\\ \noalign{\medskip}0&0&1&-84&-240
\\ \noalign{\medskip}0&0&0&20&{{254}/{7}}\\ \noalign{\medskip}0&0
&0&-21/2&-19
\end {array}\end{pmatrix},\]\[\begin{pmatrix}\begin {array}{rrrrr}2&19&-1/49&{{10}/{7}}&4
\\ \noalign{\medskip}-2/7&-{{19}/{7}}&{{1}/{343}}&-{{24}/{49}}&-{{66}/{49}}\\ \noalign{\medskip}0&0&1&-84&-240
\\ \noalign{\medskip}1&6&-{{4}/{49}}&{{26}/{7}}&{{72}/{7}}\\ \noalign{\medskip}-1/2&-3&{{2}/{49}}&-{{33}/{14}}&-6
\end{array}\end{pmatrix}.\]}

It is clear from the above computation that the $\Q$ - rank of the orthogonal group $\mathrm{SO}_\mathrm{Q}$ is two. If we denote by \[d = b a^{6} b^{-1} ,\quad
w_{_1} = b  a^{-1} b a^{-4},\quad
w_{_2} = b^{-1} a^{6} b,\quad
w_{_3} =  [w_{_2}^{-1}, w_{_1} ],\quad
w_{_4} = d w_{_1} d^{-1} w_{_1},\quad
w_{_5} = w_{_4}  w_{_3}^{-1},\]
then the final two unipotent elements, corresponding to the highest and second highest roots, are
\[q_{_1}  = a^{6},\quad
q_{_2} = w_{_5}^{7} q_{_1}^{2263},\]
which are respectively
{\tiny\[\begin{pmatrix}\begin {array}{rrrrr}1&0&0&84&0\\ \noalign{\medskip}0&1&0&0&
24\\ \noalign{\medskip}0&0&1&0&0\\ \noalign{\medskip}0&0&0&1&0
\\ \noalign{\medskip}0&0&0&0&1
\end {array}\end{pmatrix},\qquad \begin{pmatrix}\begin {array}{rrrrr}1&0&-24&0&338688\\ \noalign{\medskip}0&1
&0&0&0\\ \noalign{\medskip}0&0&1&0&-28224\\ \noalign{\medskip}0&0&0&1&0
\\ \noalign{\medskip}0&0&0&0&1
\end {array}\end{pmatrix},\]}and the arithmeticity of $\Gamma(f,g)$ follows from \cite[Theorem 3.5]{Ve87} (cf. Remarks \ref{methodoftheproof2}, \ref{methodoftheproof3}, \ref{methodoftheproof4}).

\subsection{Arithmeticity of Case \ref{arithmetic-BS17}}\label{subsection-BS17}
In this case \[\alpha=\left(0, \frac{1}{6},\frac{1}{6}, \frac{5}{6}, \frac{5}{6}\right);\qquad \beta=\left(\frac{1}{2},\frac{1}{4},\frac{1}{4},\frac{3}{4},\frac{3}{4}\right)\]
\[f(x)= x^5-3x^4+5x^3-5x^2+3x-1;\quad g(x)=x^5 +x^4+2x^3+2x^2+x+1.\]

The matrix (up to scalar multiplication) of the quadratic form $\mathrm{Q}$, 
with respect to the standard basis $\{e_1, e_2, e_3, e_4, e_5\}$, and  the change of basis matrix $X$, are
{\tiny\[\mathrm{Q}=\begin{pmatrix}\begin {array}{rrrrr}73&17&-71&-55&73\\ \noalign{\medskip}17&
73&17&-71&-55\\ \noalign{\medskip}-71&17&73&17&-71
\\ \noalign{\medskip}-55&-71&17&73&17\\ \noalign{\medskip}73&-55&-71&
17&73
\end {array}
\end{pmatrix};\quad X=\begin{pmatrix}\begin {array}{rrrrr} 2&0&0&0&1\\ \noalign{\medskip}2&0&2&-8&-
3\\ \noalign{\medskip}2&-6&2&1&6\\ \noalign{\medskip}2&3&-1&-7/2&-4
\\ \noalign{\medskip}0&-9&5&-7/2&4
\end {array}\end{pmatrix},\]}and the matrices of $\mathrm{Q}$, $A$ and $B$, with respect to the new basis, are $X^t\mathrm{Q}X$, $a=X^{-1}AX$ and $b=X^{-1}BX$; which are respectively
{\tiny\[\begin{pmatrix}\begin {array}{rrrrr} 0&0&0&0&144\\ \noalign{\medskip}0&0&0&-
432&0\\ \noalign{\medskip}0&0&144&0&0\\ \noalign{\medskip}0&-432&0&0&0
\\ \noalign{\medskip}144&0&0&0&0
\end {array}\end{pmatrix},\quad \begin{pmatrix}\begin {array}{rrrrr}0&3/2&-3/2&9/4&0\\ \noalign{\medskip}-3/
2&-6&7/2&-3/2&{{25}/{12}}\\ \noalign{\medskip}-3&-6&2&3&5/2
\\ \noalign{\medskip}-1&0&-1&3&1/2\\ \noalign{\medskip}0&-12&8&-8&4
\end {array}\end{pmatrix},\quad\begin{pmatrix}\begin {array}{rrrrr}0&6&-4&4&-2\\ \noalign{\medskip}-3/2&-6&
7/2&-3/2&{{25}/{12}}\\ \noalign{\medskip}-3&-6&2&3&5/2
\\ \noalign{\medskip}-1&0&-1&3&1/2\\ \noalign{\medskip}0&-3&3&-9/2&0
\end{array}\end{pmatrix}.\]}

It is clear from the above computation that the $\Q$ - rank of the orthogonal group $\mathrm{SO}_\mathrm{Q}$ is two. If we denote by \[d = b a^{-6} b^{-1},\quad
e = a^{-1} b^{-4} a,\quad
w_{_1} = (a  b^{-1} a)^{8},\quad w_{_2} = w_{1}^{2} e,\quad w_{_3} = w_{_2}^{-3} w_{_1}^{2},\quad
w_{_4} = d w_{_2} d^{-1},\]
then the final three unipotent elements are
\[q_{_1}  = [ w_{_4}, w_{_2}^{-7}],\quad
q_{_2} =w_{_2}^{4032} q_{_1}^{30},\quad
q_{_3}  = w_{_3}^{112} q_{_1}^{-1},\]
which are respectively
{\tiny\[\begin{pmatrix}\begin {array}{rrrrr}1&0&0&-1344&0\\ \noalign{\medskip}0&1&0&0
&-448\\ \noalign{\medskip}0&0&1&0&0\\ \noalign{\medskip}0&0&0&1&0
\\ \noalign{\medskip}0&0&0&0&1
\end {array}\end{pmatrix},\quad \begin{pmatrix}\begin {array}{rrrrr}1&0&-16128&0&-130056192
\\ \noalign{\medskip}0&1&0&0&0\\ \noalign{\medskip}0&0&1&0&16128
\\ \noalign{\medskip}0&0&0&1&0\\ \noalign{\medskip}0&0&0&0&1
\end {array}\end{pmatrix},\quad \begin{pmatrix}\begin {array}{rrrrr}1&1344&0&0&0\\ \noalign{\medskip}0&1&0&0
&0\\ \noalign{\medskip}0&0&1&0&0\\ \noalign{\medskip}0&0&0&1&448
\\ \noalign{\medskip}0&0&0&0&1
\end {array}\end{pmatrix}.\]}

\subsection{Arithmeticity of Case \ref{arithmetic-BS18}}\label{subsection-BS18}
In this case \[\alpha=\left(0, \frac{1}{8},\frac{3}{8}, \frac{5}{8}, \frac{7}{8}\right);\qquad \beta=\left(\frac{1}{2},\frac{1}{3},\frac{1}{3},\frac{2}{3},\frac{2}{3}\right)\]
\[f(x)= x^5-x^4+x-1;\quad g(x)=x^5 +3x^4+5x^3+5x^2+3x+1.\]

The matrix (up to scalar multiplication) of the quadratic form $\mathrm{Q}$, 
with respect to the standard basis $\{e_1, e_2, e_3, e_4, e_5\}$, and  the change of basis matrix $X$, are
{\tiny\[\mathrm{Q}=\begin{pmatrix}\begin {array}{rrrrr}73&-53&1&55&-71\\ \noalign{\medskip}-53&
73&-53&1&55\\ \noalign{\medskip}1&-53&73&-53&1\\ \noalign{\medskip}55&
1&-53&73&-53\\ \noalign{\medskip}-71&55&1&-53&73
\end {array}
\end{pmatrix};\quad X=\begin{pmatrix}\begin {array}{rrrrr}-2&{{71}/{7}}&-{{4}/{49}}&9&7
\\ \noalign{\medskip}0&-{{50}/{7}}&{{4}/{49}}&-2&-5
\\ \noalign{\medskip}1&-3&0&-7&-1\\ \noalign{\medskip}2&{{26}/{7}
}&{{8}/{49}}&-11&4\\ \noalign{\medskip}-1&{{118}/{7}}&{{4}/{49}}&5&13
\end {array}\end{pmatrix},\]}and the matrices of $\mathrm{Q}$, $A$ and $B$, with respect to the new basis, are $X^t\mathrm{Q}X$, $a=X^{-1}AX$ and $b=X^{-1}BX$; which are respectively
{\tiny\[\begin{pmatrix}\begin {array}{rrrrr}0&0&0&0&-36\\ \noalign{\medskip}0&0&0&-
252&0\\ \noalign{\medskip}0&0&{{144}/{49}}&0&0
\\ \noalign{\medskip}0&-252&0&0&0\\ \noalign{\medskip}-36&0&0&0&0
\end {array}\end{pmatrix},\qquad \begin{pmatrix}\begin {array}{rrrrr}3&-{{24}/{7}}&{{12}/{49}}&-13&
-2\\ \noalign{\medskip}-5/7&{{64}/{49}}&-{{60}/{343}}&{{30}/{7}}&6/7\\ \noalign{\medskip}1&{{46}/{7}}&-{{37}/{49}}&-6&3\\ \noalign{\medskip}{{29}/{49}}&{{11}/{343}}&{{152}/{2401}}&-{{125}/{49}}&{{3}/{49}}
\\ \noalign{\medskip}1&-3/7&{{12}/{49}}&-6&0
\end {array}\end{pmatrix},\]\[\begin{pmatrix}\begin {array}{rrrrr}-4&{{802}/{7}}&{{40}/{49}}&22&
89\\ \noalign{\medskip}-2&{{1126}/{49}}&-{{24}/{343}}&{{75}/{7}}&{{123}/{7}}\\ \noalign{\medskip}0&{{164}/{7}}
&-{{33}/{49}}&-1&16\\ \noalign{\medskip}-6/7&{{8389}/{343}}&
{{436}/{2401}}&{{230}/{49}}&{{926}/{49}}
\\ \noalign{\medskip}3&-{{239}/{7}}&{{4}/{49}}&-16&-26
\end{array}\end{pmatrix}.\]}

It is clear from the above computation that the $\Q$ - rank of the orthogonal group $\mathrm{SO}_\mathrm{Q}$ is two. If we denote by \[c = a^{-1} b,\ 
d = c b c,\quad
e = a^{-1} b^{6} a,\quad
w_{_1} = [ e, a^{4}],\quad
w_{_2} = w_{1}^{49} e^{-48},\quad
w_{_3} = e d^{-3},\ 
w_{_4} = w_{_3} w_{_1}^{-1}  w_{_3} w_{_1},\]
then the final three unipotent elements are
\[q_{_1}  =q_{_2}^{-1} w_{_2}^{6},\quad
q_{_2} =  w_{_4}^{7} e^{-247},\quad
q_{_3}  = e,\]
which are respectively
{\tiny\[\begin{pmatrix}\begin {array}{rrrrr}1&0&0&-3528&0\\ \noalign{\medskip}0&1&0&0
&504\\ \noalign{\medskip}0&0&1&0&0\\ \noalign{\medskip}0&0&0&1&0
\\ \noalign{\medskip}0&0&0&0&1
\end {array}\end{pmatrix},\quad \begin{pmatrix}\begin {array}{rrrrr}1&0&144&0&127008\\ \noalign{\medskip}0&1
&0&0&0\\ \noalign{\medskip}0&0&1&0&1764\\ \noalign{\medskip}0&0&0&1&0
\\ \noalign{\medskip}0&0&0&0&1
\end {array}\end{pmatrix},\quad \begin{pmatrix}\begin {array}{rrrrr}1&42&0&0&0\\ \noalign{\medskip}0&1&0&0&0
\\ \noalign{\medskip}0&0&1&0&0\\ \noalign{\medskip}0&0&0&1&-6
\\ \noalign{\medskip}0&0&0&0&1
\end {array}\end{pmatrix}.\]}

\subsection{Arithmeticity of Case \ref{arithmetic-BS19}}\label{subsection-BS19}
In this case \[\alpha=\left(0, \frac{1}{8},\frac{3}{8}, \frac{5}{8}, \frac{7}{8}\right);\qquad \beta=\left(\frac{1}{2},\frac{1}{3},\frac{2}{3},\frac{1}{4},\frac{3}{4}\right)\]
\[f(x)=x^5-x^4+x-1;\quad g(x)=x^5 +2x^4+3x^3+3x^2+2x+1.\]

The matrix (up to scalar multiplication) of the quadratic form $\mathrm{Q}$, 
with respect to the standard basis $\{e_1, e_2, e_3, e_4, e_5\}$, and  the change of basis matrix $X$, are
{\tiny\[\mathrm{Q}=\begin{pmatrix}\begin {array}{rrrrr}13&-11&1&13&-11\\ \noalign{\medskip}-11&
13&-11&1&13\\ \noalign{\medskip}1&-11&13&-11&1\\ \noalign{\medskip}13&
1&-11&13&-11\\ \noalign{\medskip}-11&13&1&-11&13
\end {array}
\end{pmatrix};\quad X=\begin{pmatrix}\begin {array}{rrrrr}-8&4&-3/2&-17/2&1\\ \noalign{\medskip}-8
&0&-3&-12&1\\ \noalign{\medskip}-8&0&-3&-16&1\\ \noalign{\medskip}2&-4
&3/2&1/2&-2\\ \noalign{\medskip}-2&0&0&0&-1
\end {array}\end{pmatrix},\]}and the matrices of $\mathrm{Q}$, $A$ and $B$, with respect to the new basis, are $X^t\mathrm{Q}X$, $a=X^{-1}AX$ and $b=X^{-1}BX$; which are respectively
{\tiny\[\begin{pmatrix}\begin {array}{rrrrr} 0&0&0&0&48\\ \noalign{\medskip}0&0&0&-
192&0\\ \noalign{\medskip}0&0&36&0&0\\ \noalign{\medskip}0&-192&0&0&0
\\ \noalign{\medskip}48&0&0&0&0
\end {array}\end{pmatrix},\qquad \begin{pmatrix}\begin {array}{rrrrr}3/2&-3&3/8&{{19}/{8}}&-5/4
\\ \noalign{\medskip}{{39}/{16}}&-{{33}/{8}}&{{69}/{64}
}&{{297}/{64}}&-{{73}/{32}}\\ \noalign{\medskip}-5&6&-11/4&-
{{35}/{4}}&7/2\\ \noalign{\medskip}1/2&1&3/8&{{7}/{8}}&1/4
\\ \noalign{\medskip}-3&10&-9/4&-{{21}/{4}}&11/2
\end {array}\end{pmatrix},\] \[\begin{pmatrix}\begin {array}{rrrrr}5/2&-3&3/8&{{19}/{8}}&-3/4
\\ \noalign{\medskip}{{73}/{16}}&-{{33}/{8}}&{{69}/{64}
}&{{297}/{64}}&-{{39}/{32}}\\ \noalign{\medskip}-7&6&-11/4&-
{{35}/{4}}&5/2\\ \noalign{\medskip}-1/2&1&3/8&{{7}/{8}}&-1/4
\\ \noalign{\medskip}-11&10&-9/4&-{{21}/{4}}&3/2
\end{array}\end{pmatrix}.\]}

It is clear from the above computation that the $\Q$ - rank of the orthogonal group $\mathrm{SO}_\mathrm{Q}$ is two. If we denote by \[c = a^{-1} b,\quad
d = [a^{-4}, b^{3}],\quad
e = a b,\quad
w_{_1} = c d c,\quad
w_{_2} = c e,\quad
w_{_3} = w_{_2} w_{_1} w_{_2}^{-1} ,\quad
w_{_4} = w_{_1}^{2}  w_{_3},\]
then the final two unipotent elements, corresponding to the highest and second highest roots, are
\[q_{_1}  = [ w_{_1} , w_{_3}],\quad
q_{_2} = w_{_4}^{-96} q_{_1}^{132},\]
which are respectively
{\tiny\[\begin{pmatrix}\begin {array}{rrrrr}1&0&0&96&0\\ \noalign{\medskip}0&1&0&0&
24\\ \noalign{\medskip}0&0&1&0&0\\ \noalign{\medskip}0&0&0&1&0
\\ \noalign{\medskip}0&0&0&0&1
\end {array}\end{pmatrix},\qquad \begin{pmatrix}\begin {array}{rrrrr}1&0&-576&0&-221184\\ \noalign{\medskip}0
&1&0&0&0\\ \noalign{\medskip}0&0&1&0&768\\ \noalign{\medskip}0&0&0&1&0
\\ \noalign{\medskip}0&0&0&0&1
\end {array}\end{pmatrix},\]}and the arithmeticity of $\Gamma(f,g)$ follows from \cite[Theorem 3.5]{Ve87} (cf. Remarks \ref{methodoftheproof2}, \ref{methodoftheproof3}, \ref{methodoftheproof4}).

\subsection{Arithmeticity of Case \ref{arithmetic-BS20}}\label{subsection-BS20}
In this case \[\alpha=\left(0, \frac{1}{8},\frac{3}{8}, \frac{5}{8}, \frac{7}{8}\right);\qquad \beta=\left(\frac{1}{2},\frac{1}{4},\frac{1}{4},\frac{3}{4},\frac{3}{4}\right)\]
\[f(x)=x^5-x^4+x-1;\quad g(x)=x^5 +x^4+2x^3+2x^2+x+1.\]

The matrix (up to scalar multiplication) of the quadratic form $\mathrm{Q}$, 
with respect to the standard basis $\{e_1, e_2, e_3, e_4, e_5\}$, and  the change of basis matrix $X$, are
{\tiny\[\mathrm{Q}=\begin{pmatrix}\begin {array}{rrrrr}1&-3&1&5&1\\ \noalign{\medskip}-3&1&-3&1
&5\\ \noalign{\medskip}1&-3&1&-3&1\\ \noalign{\medskip}5&1&-3&1&-3
\\ \noalign{\medskip}1&5&1&-3&1
\end {array}
\end{pmatrix};\quad X=\begin{pmatrix}\begin {array}{rrrrr}0&1&{{8}/{9}}&1/6&1
\\ \noalign{\medskip}1&1&{{8}/{9}}&1/6&1/2\\ \noalign{\medskip}0&
4&{{8}/{9}}&-1/3&2\\ \noalign{\medskip}-1&1&{{8}/{9}}&-{{17}/{6}}&1/2\\ \noalign{\medskip}0&1&{{8}/{9}}&1/6&0
\end {array}\end{pmatrix},\]}and the matrices of $\mathrm{Q}$, $A$ and $B$, with respect to the new basis, are $X^t\mathrm{Q}X$, $a=X^{-1}AX$ and $b=X^{-1}BX$; which are respectively
{\tiny\[\begin{pmatrix}\begin {array}{rrrrr}0&0&0&0&-8\\ \noalign{\medskip}0&0&0&24&0
\\ \noalign{\medskip}0&0&{{64}/{9}}&0&0\\ \noalign{\medskip}0&24&0
&0&0\\ \noalign{\medskip}-8&0&0&0&0
\end {array}\end{pmatrix}, \begin{pmatrix}\begin {array}{rrrrr}1/2&-3/2&-4/3&5/4&3/4
\\ \noalign{\medskip}-1/18&{{5}/{18}}&{{32}/{81}}&-{{115}/{108}}&{{7}/{36}}\\ \noalign{\medskip}-1&2&{{13}/{9}}&-5
/3&1/2\\ \noalign{\medskip}-1/3&-1/3&{{16}/{27}}&-{{13}/{18}
}&-5/6\\ \noalign{\medskip}1&-1&-{{8}/{9}}&{{17}/{6}}&-1/2
\end {array}\end{pmatrix},\begin{pmatrix}\begin {array}{rrrrr}1/2&1/2&4/9&{{19}/{12}}&3/4
\\ \noalign{\medskip}-1/18&1/6&{{8}/{27}}&-{{13}/{12}}&{{7}/{36}}\\ \noalign{\medskip}-1&0&-1/3&-2&1/2
\\ \noalign{\medskip}-1/3&-1&0&-5/6&-5/6\\ \noalign{\medskip}1&-1&-{{8}/{9}}&{{17}/{6}}&-1/2
\end{array}\end{pmatrix}.\]}

It is clear from the above computation that the $\Q$ - rank of the orthogonal group $\mathrm{SO}_\mathrm{Q}$ is two. If we denote by \[w_{_1} = [ a,  b^{-1}],\quad
 w_{_2} = [ a^{-4}, b^{3}],\quad
 w_{_3} =  b a^{4} b^{-1},\quad 
 w_{_4} = [w_{_3}, w_{_1}],\quad
 w_{_5} = w_{_2} w_{_1} w_{_2}^{-1} w_{_1}^{-89},\]
then the final three unipotent elements are
\[q_{_1}=w_{5}^{3} q_{_2}^{-8},\quad
q_{_2} =  w_{5} w_{4}^{-10},\quad
q_{_3}=(w_{1}^{18} q_{_2})^{40} q_{_1},\]
which are respectively
{\tiny\[\begin{pmatrix}\begin {array}{rrrrr}1&0&0&240&0\\ \noalign{\medskip}0&1&0&0&
80\\ \noalign{\medskip}0&0&1&0&0\\ \noalign{\medskip}0&0&0&1&0
\\ \noalign{\medskip}0&0&0&0&1
\end {array}\end{pmatrix},\quad \begin{pmatrix}\begin {array}{rrrrr}1&0&32&0&576\\ \noalign{\medskip}0&1&0&0
&0\\ \noalign{\medskip}0&0&1&0&36\\ \noalign{\medskip}0&0&0&1&0
\\ \noalign{\medskip}0&0&0&0&1
\end {array}\end{pmatrix},\quad \begin{pmatrix}\begin {array}{rrrrr}1&-1440&0&0&0\\ \noalign{\medskip}0&1&0&0
&0\\ \noalign{\medskip}0&0&1&0&0\\ \noalign{\medskip}0&0&0&1&-480
\\ \noalign{\medskip}0&0&0&0&1
\end {array}\end{pmatrix}.\]}

\subsection{Arithmeticity of Case \ref{arithmetic-BS21} ($\Q$ - rank one case)}\label{subsection-BS21}
In this case \[\alpha=\left(0, \frac{1}{12},\frac{5}{12}, \frac{7}{12}, \frac{11}{12}\right);\qquad \beta=\left(\frac{1}{2},\frac{1}{3},\frac{2}{3},\frac{1}{4},\frac{3}{4}\right)\]
\[f(x)=x^5-x^4-x^3+x^2+x-1;\quad g(x)=x^5 +2x^4+3x^3+3x^2+2x+1.\]

The matrix (up to scalar multiplication) of the quadratic form $\mathrm{Q}$, 
with respect to the standard basis $\{e_1, e_2, e_3, e_4, e_5\}$, and  the change of basis matrix $X$, are
{\tiny\[\mathrm{Q}=\begin{pmatrix}\begin {array}{rrrrr}5&-5&3&1&-1\\ \noalign{\medskip}-5&5&-5&
3&1\\ \noalign{\medskip}3&-5&5&-5&3\\ \noalign{\medskip}1&3&-5&5&-5
\\ \noalign{\medskip}-1&1&3&-5&5
\end {array}
\end{pmatrix};\quad X=\begin{pmatrix}\begin {array}{rrrrr}0&-3&5&2&2\\ \noalign{\medskip}0&-4&9&1&
1\\ \noalign{\medskip}1&-5&13&0&1\\ \noalign{\medskip}1&-4&10&2&-3
\\ \noalign{\medskip}0&-2&1&3&-3
\end {array}\end{pmatrix},\]}and the matrices of $\mathrm{Q}$, $A$ and $B$, with respect to the new basis, are $X^t\mathrm{Q}X$, $a=X^{-1}AX$ and $b=X^{-1}BX$; which are respectively
{\tiny\[\begin{pmatrix}\begin {array}{rrrrr} 0&0&0&0&12\\ \noalign{\medskip}0&24&-12&
-36&0\\ \noalign{\medskip}0&-12&-24&60&0\\ \noalign{\medskip}0&-36&60&
24&0\\ \noalign{\medskip}12&0&0&0&0
\end {array}\end{pmatrix}, \begin{pmatrix}\begin {array}{rrrrr}1/6&-5/2&{{19}/{6}}&8/3&-{{49}/{6}}\\ \noalign{\medskip}0&0&-4&1&0\\ \noalign{\medskip}0&0&-1&0&1
\\ \noalign{\medskip}1/6&-3/2&-5/6&8/3&-{{19}/{6}}
\\ \noalign{\medskip}-1/6&1/2&-{{13}/{6}}&1/3&-5/6
\end {array}\end{pmatrix}, \begin{pmatrix}\begin {array}{rrrrr}1/6&3/2&7/6&-10/3&-{{13}/{6}}
\\ \noalign{\medskip}0&0&-4&1&0\\ \noalign{\medskip}0&0&-1&0&1
\\ \noalign{\medskip}1/6&1/2&-{{11}/{6}}&-1/3&-1/6
\\ \noalign{\medskip}-1/6&1/2&-{{13}/{6}}&1/3&-5/6
\end{array}\end{pmatrix}.\]}

A computation shows that the $\Q$ - rank of the orthogonal group $\mathrm{SO}_\mathrm{Q}$ is one. If we denote by \[c = a^{-1} b,\quad
d = [a^{-1}, b^{-5}],\quad
e =  b a b^{-1},\quad
r = [ b^{-1} , a],\quad
s = b^{-1} a b ,\quad
w_{_1} = [a^{-1}, b^{2}],\]
\[w_{_2} = c [d, w_{_1}] c,\quad
w_{_3} = [ e, r^{-1}] ,\quad
w_{_4} = w_{_2} w_{_3} w_{_2},\quad
w_{_5} = [ s, r^{-1} ],\quad
w_{_6} = w_{_4}  w_{_5}^{-1} w_{_4}  w_{_5}^{4},\]
\[w_{_7} = [ w_{_4}, r^{-1} ] w_{_5}^{39},\quad
w_{_8} = [ w_{_1}, r^{-1} ],\]
then the final three unipotent elements are
\[q_{_1}  = w_{_6}^{504} w_{_7}^{306},\quad
q_{_2} = w_{_6}^{360} w_{_7}^{342},\quad
q_{_3}  = w_{_8}^{62208} q_{_2}^{-78},\]which are respectively
{\tiny\[\begin{pmatrix}\begin {array}{rrrrr}1&0&0&62208&-806215680
\\ \noalign{\medskip}0&1&0&0&-57024\\ \noalign{\medskip}0&0&1&0&-36288
\\ \noalign{\medskip}0&0&0&1&-25920\\ \noalign{\medskip}0&0&0&0&1
\end {array}\end{pmatrix}, \begin{pmatrix}\begin {array}{rrrrr}1&0&62208&0&-806215680
\\ \noalign{\medskip}0&1&0&0&-67392\\ \noalign{\medskip}0&0&1&0&-25920
\\ \noalign{\medskip}0&0&0&1&-36288\\ \noalign{\medskip}0&0&0&0&1
\end {array}\end{pmatrix}, \begin{pmatrix}\begin {array}{rrrrr}1&-1866240&0&0&-4208445849600
\\ \noalign{\medskip}0&1&0&0&4510080\\ \noalign{\medskip}0&0&1&0&
2021760\\ \noalign{\medskip}0&0&0&1&1710720\\ \noalign{\medskip}0&0&0&0
&1
\end {array}\end{pmatrix}.\]}

\subsection{Arithmeticity of Case \ref{arithmetic-BS22}}\label{subsection-BS22}
In this case \[\alpha=\left(0, \frac{1}{12},\frac{5}{12}, \frac{7}{12}, \frac{11}{12}\right);\qquad \beta=\left(\frac{1}{2},\frac{1}{5},\frac{2}{5},\frac{3}{5},\frac{4}{5}\right)\]
\[f(x)=x^5-x^4-x^3+x^2+x-1;\quad g(x)=x^5+2x^4+2x^3+2x^2+2x+1.\]

The matrix (up to scalar multiplication) of the quadratic form $\mathrm{Q}$, 
with respect to the standard basis $\{e_1, e_2, e_3, e_4, e_5\}$, and  the change of basis matrix $X$, are
{\tiny\[\mathrm{Q}=\begin{pmatrix}\begin {array}{rrrrr}21&-19&11&1&-9\\ \noalign{\medskip}-19&
21&-19&11&1\\ \noalign{\medskip}11&-19&21&-19&11\\ \noalign{\medskip}1
&11&-19&21&-19\\ \noalign{\medskip}-9&1&11&-19&21
\end {array}
\end{pmatrix};\quad X=\begin{pmatrix}\begin {array}{rrrrr}-3&-2&-5/4&7/4&-{{15}/{4}}
\\ \noalign{\medskip}-4&-1&0&3&-4\\ \noalign{\medskip}-2&1&5/2&7/2&-3/
2\\ \noalign{\medskip}5&2&{{15}/{4}}&7/4&1/4\\ \noalign{\medskip}
4&0&0&0&-1
\end {array}\end{pmatrix},\]}and the matrices of $\mathrm{Q}$, $A$ and $B$, with respect to the new basis, are $X^t\mathrm{Q}X$, $a=X^{-1}AX$ and $b=X^{-1}BX$; which are respectively
{\tiny\[\begin{pmatrix}\begin {array}{rrrrr} 0&0&0&0&-40\\ \noalign{\medskip}0&0&0&-
10&0\\ \noalign{\medskip}0&0&25&0&0\\ \noalign{\medskip}0&-10&0&0&0
\\ \noalign{\medskip}-40&0&0&0&0
\end {array}\end{pmatrix},\qquad \begin{pmatrix}\begin {array}{rrrrr}{{33}/{8}}&{{17}/{16}}&{{45}/{32}}&{{9}/{32}}&{{3}/{32}}\\ \noalign{\medskip}-{{25}/{2}}&-9/4&-5/8&{{15}/{8}}&-{{11}/{8}}
\\ \noalign{\medskip}-3&-3/2&-11/4&-3/4&-1/4\\ \noalign{\medskip}9&3&{{15}/{4}}&3/4&1/4\\ \noalign{\medskip}15/2&9/4&{{15}/{8}}&-5
/8&{{9}/{8}}
\end {array}\end{pmatrix},\]\[\begin{pmatrix}\begin {array}{rrrrr}3/8&{{17}/{16}}&{{45}/{32}}&{{9}/{32}}&{{33}/{32}}\\ \noalign{\medskip}-11/2&-9/4&-5/8&{{15}/{8}}&-{{25}/{8}}\\ \noalign{\medskip}-1&-3/2&-11/4&-3/4
&-3/4\\ \noalign{\medskip}1&3&{{15}/{4}}&3/4&9/4
\\ \noalign{\medskip}9/2&9/4&{{15}/{8}}&-5/8&{{15}/{8}}
\end{array}\end{pmatrix}.\]}

It is clear from the above computation that the $\Q$ - rank of the orthogonal group $\mathrm{SO}_\mathrm{Q}$ is two. If we denote by \[c = a^{-1} b,\quad
d = b^{-4} a^{6} b^{4},\quad
e =  d a^{4} d a^{-4},\quad
w_{_1} = [e, a^{6}],\quad
w_{_2} = [w_{_1}, a^{3}],\quad
w_{_3} = [ w_{_1}^{-24}, w_{_2}],\]
\[w_{_4} = w_{_1}^{-384} w_{_3},\ 
w_{_5} = w_{_2}^{96} w_{_4},\quad
w_{_6} = c w_{_3} c,\ 
w_{_7} = d a^{6} d^{-1} a^6,\quad
w_{_8} = w_{_6}^{5} w_{_7} w_{_3}^{-1}w_{_7}^{-1},\ 
w_{_9} = w_{_8}^{-1} w_{_5} w_{_8},\]
then the final two unipotent elements, corresponding to the highest and second highest roots, are
\[q_{_1}  = w_{_3},\quad
q_{_2} = w_{_9}^{3840} q_{_1}^{238886400},\]
which are respectively
{\tiny\[\begin{pmatrix}\begin {array}{rrrrr}1&0&0&-3840&0\\ \noalign{\medskip}0&1&0&0
&15360\\ \noalign{\medskip}0&0&1&0&0\\ \noalign{\medskip}0&0&0&1&0
\\ \noalign{\medskip}0&0&0&0&1
\end {array}\end{pmatrix},\qquad \begin{pmatrix}\begin {array}{rrrrr}1&0&7372800&0&43486543872000
\\ \noalign{\medskip}0&1&0&0&0\\ \noalign{\medskip}0&0&1&0&11796480
\\ \noalign{\medskip}0&0&0&1&0\\ \noalign{\medskip}0&0&0&0&1
\end {array}\end{pmatrix},\]}and the arithmeticity of $\Gamma(f,g)$ follows from \cite[Theorem 3.5]{Ve87} (cf. Remarks \ref{methodoftheproof2}, \ref{methodoftheproof3}, \ref{methodoftheproof4}).

\subsection{Arithmeticity of Case \ref{arithmetic-BS23}}\label{subsection-BS23}
In this case \[\alpha=\left(0, \frac{1}{12},\frac{5}{12}, \frac{7}{12}, \frac{11}{12}\right);\qquad \beta=\left(\frac{1}{2},\frac{1}{8},\frac{3}{8},\frac{5}{8},\frac{7}{8}\right)\]
\[f(x)=x^5-x^4-x^3+x^2+x-1;\quad g(x)=x^5+x^4+x+1.\]

The matrix (up to scalar multiplication) of the quadratic form $\mathrm{Q}$, 
with respect to the standard basis $\{e_1, e_2, e_3, e_4, e_5\}$, and  the change of basis matrix $X$, are
{\tiny\[\mathrm{Q}=\begin{pmatrix}\begin {array}{rrrrr}1&-3&1&1&1\\ \noalign{\medskip}-3&1&-3&1
&1\\ \noalign{\medskip}1&-3&1&-3&1\\ \noalign{\medskip}1&1&-3&1&-3
\\ \noalign{\medskip}1&1&1&-3&1
\end {array}
\end{pmatrix};\quad X=\begin{pmatrix}\begin {array}{rrrrr} 0&1&3/4&-1/6&2\\ \noalign{\medskip}1&-1&
-3/4&1/6&-1\\ \noalign{\medskip}0&-1&-3/4&3/2&-2\\ \noalign{\medskip}-
1&0&1&0&2\\ \noalign{\medskip}0&1&3/4&-1/6&3
\end {array}\end{pmatrix},\]}and the matrices of $\mathrm{Q}$, $A$ and $B$, with respect to the new basis, are $X^t\mathrm{Q}X$, $a=X^{-1}AX$ and $b=X^{-1}BX$; which are respectively
{\tiny\[\begin{pmatrix}\begin {array}{rrrrr}0&0&0&0&4\\ \noalign{\medskip}0&0&0&16/3
&0\\ \noalign{\medskip}0&0&1&0&0\\ \noalign{\medskip}0&16/3&0&0&0
\\ \noalign{\medskip}4&0&0&0&0
\end {array}\end{pmatrix},\quad\begin{pmatrix}\begin {array}{rrrrr}1&1&-1/4&-1/6&0\\ \noalign{\medskip}-1/8
&1/8&{{11}/{32}}&-{{49}/{48}}&{{9}/{8}}
\\ \noalign{\medskip}3&1&-9/4&7/6&-3\\ \noalign{\medskip}3/4&-3/4&-{{9}/{16}}&1/8&-3/4\\ \noalign{\medskip}-1&0&1&0&2
\end {array}\end{pmatrix},\quad\begin{pmatrix}\begin {array}{rrrrr}1&-1&-7/4&1/6&-6\\ \noalign{\medskip}-1/
8&1/4&{{7}/{16}}&-{{25}/{24}}&3/2\\ \noalign{\medskip}3&-2&-
9/2&5/3&-12\\ \noalign{\medskip}3/4&-3/2&-{{9}/{8}}&1/4&-3
\\ \noalign{\medskip}-1&0&1&0&2
\end{array}\end{pmatrix}.\]}

It is clear from the above computation that the $\Q$ - rank of the orthogonal group $\mathrm{SO}_\mathrm{Q}$ is two. If we denote by \[w_{_1} = [a^{-1}, b],\quad
w_{_2} = a^3 b^{4} a^{-3},\quad
w_{_3} = (a^3 b^{4})^{4},\quad
w_{_4} = w_{_2} w_{_1} w_{_2},\]
then the final three unipotent elements are
\[q_{_1}=q_{_2}^{2} w_{_4}^{24} w_{_3}^{-6},\quad
q_{_2} =w_{_1}^{12} w_{_3}^{-3}  ,\quad
q_{_3}=(w_{_3}^{12} q_{_2})^{-8} q_{_1},\]
which are respectively
{\tiny\[\begin{pmatrix}\begin {array}{rrrrr}1&0&0&-64&0\\ \noalign{\medskip}0&1&0&0&
48\\ \noalign{\medskip}0&0&1&0&0\\ \noalign{\medskip}0&0&0&1&0
\\ \noalign{\medskip}0&0&0&0&1
\end {array}\end{pmatrix},\quad \begin{pmatrix}\begin {array}{rrrrr}1&0&12&0&-288\\ \noalign{\medskip}0&1&0&0
&0\\ \noalign{\medskip}0&0&1&0&-48\\ \noalign{\medskip}0&0&0&1&0
\\ \noalign{\medskip}0&0&0&0&1
\end {array}\end{pmatrix},\quad \begin{pmatrix}\begin {array}{rrrrr}1&-384&0&0&0\\ \noalign{\medskip}0&1&0&0
&0\\ \noalign{\medskip}0&0&1&0&0\\ \noalign{\medskip}0&0&0&1&288
\\ \noalign{\medskip}0&0&0&0&1
\end {array}\end{pmatrix}.\]}

%%%%%%%%%%%%%%%%%%%%%%%%%%%%%%%%%%%%%%%%%%%%%%%%%%
%%% BIBLIOGRAPHY

\end{document}